\documentclass[reqno,11pt]{amsart}

\newtheorem{theorem}{Theorem}[section]
\newtheorem{lemma}[theorem]{Lemma}
\newtheorem{corollary}[theorem]{Corollary}

\theoremstyle{definition}

\newtheorem{assumption}{Assumption}[section]

 \theoremstyle{remark}
\newtheorem{remark}[theorem]{Remark}

\newcommand\bL{\mathbb{L}}

\newcommand\bR{\mathbb{R}}

\newcommand\bW{\mathbb{W}}

\newcommand\cB{\mathcal{B}}

\newcommand\cF{\mathcal{F}}
\newcommand\cH{\mathcal{H}}
 
\newcommand\cP{\mathcal{P}}

\newcommand\cW{\mathcal{W}}

\newcommand\frW{\mathfrak{W}}

\newcommand{\mysection}[1]{\section{#1}
\setcounter{equation}{0}}

\newcommand\cbrk{\text{$]$\kern-.15em$]$}} 
\newcommand\opar{
\text{\,\raise.2ex\hbox{${\scriptstyle |}$}\kern-.34em$($}} 
\newcommand\cpar{%
\text{$)$\kern-.34em\raise.2ex\hbox{${\scriptstyle |}$}}\,}
\newcommand\obrk{\text{$[$\kern-.15em$[$}}

\begin{document}

\title[SPDEs with VMO coefficients]
{On divergence
form SPDEs
with VMO coefficients in a half space}
\author{N.V. Krylov}
\address{127 Vincent Hall, University of Minnesota, Minneapolis,
 MN, 55455}
\thanks{The work   was partially supported by
NSF Grant DMS-0653121}
\email{krylov@math.umn.edu}

\renewcommand{\subjclassname}{%
\textup{2000} Mathematics Subject Classification}


\begin{abstract}
We extend several known results on  solvability in the 
Sobolev spaces
$W^{1}_{p}$, $p\in[2,\infty)$, of SPDEs in divergence
form in $\bR^{d}_{+}$ to equations having
coefficients which are discontinuous in the space variable.
\end{abstract}

\maketitle

\mysection{Introduction}

Let $(\Omega,\cF,P)$ be a complete probability space
with an increasing filtration $\{\cF_{t},t\geq0\}$
of complete with respect to $(\cF,P)$ $\sigma$-fields
$\cF_{t}\subset\cF$. Denote by $\cP$ the predictable
$\sigma$-field in $\Omega\times(0,\infty)$
associated with $\{\cF_{t}\}$. Let
 $w^{k}_{t}$, $k=1,2,...$, be independent one-dimensional
Wiener processes with respect to $\{\cF_{t}\}$.  

We fix a stopping time $\tau$ and for $t\leq \tau $
in  
$$
\bR^{d}_{+}=\{x=(x^{1},x'):x^{1}>0,x'=(x^{2},...,x^{d})
\in\bR^{d-1}\},\quad d\geq2,
$$
$$
\bR^{1}_{+}=
\bR_{+}=(0,\infty)
$$
 consider  the following
equation
\begin{equation}
                                          \label{11.13.1}
du_{t}=(L_{t}u_{t} 
+D_{i}f^{i}_{t}+f^{0}_{t})\,dt
+(\Lambda^{k}_{t}u_{t}+g^{k}_{t})\,dw^{k}_{t},
\end{equation}
  where $u_{t}=u_{t}(x)=u_{t}(\omega,x)$ is an unknown function,
$$
L_{t}\psi(x)=D_{j}\big(a^{ij}_{t}(x)D_{i}\psi(x)
+a^{j}_{t}(x)\psi(x)
\big)+b^{i}_{t}(x)D_{i}\psi(x)+c_{t}(x)\psi(x),
$$
$$
\Lambda^{k}_{t}\psi(x)=\sigma^{ik}_{t}(x)D_{i}\psi(x)
+\nu^{k}_{t}(x)\psi(x),
$$
the summation  convention
with respect to $i,j=1,...,d$ and $k=1,2,...$ is enforced
and detailed assumptions on the coefficients and the 
free terms will be given later. Equation \eqref{11.13.1}
is supplemented with zero initial data and {\em zero\/}
boundary condition on $x^{1}=0$. Other initial conditions
can also be considered by a standard method of continuing
them for $t>0$ and subtracting the result of
continuation from $u$. However, for simplicity
of presentation we confine ourselves to the simplest
case of zero initial condition.

One of possible approaches
to equation \eqref{11.13.1} is to rewrite it in
the  nondivergence form
assuming that the coefficients $a^{ij}_{t}$ and 
$a^{i}_{t}$ are differentiable
 in $x$ and then one could apply the results
from \cite{Ki2} to obtain the solvability in
$W^{1}_{p}$ spaces for all $p\geq2$.
It turns out that the 
differentiability of $a^{ij}_{t}$ and 
$a^{i}_{t}$ is not needed
for the corresponding 
counterparts of the results in \cite{Ki2} to be true,
which is shown in \cite{Ki}, where the coefficients
$a$ and $\sigma$ are just continuous in $x$.
Recent development in the theory of parabolic
PDEs allows one to further reduce the regularity
assumption on $a$ (but not $\sigma$) and require
that $a$ be in VMO with respect to the space variable
and showing this is the main purpose of this article.

The main guidelines we follow are quite common:
getting a priori estimates and using the method of
continuity. The method of continuity requires
a starting point, which in our case is 
the {\em solvability\/} of
the equation
\begin{equation}
                                                     \label{6.4.1}
du_{t}=(\Delta u_{t}+D_{i}f^{i}_{t}+f^{0}_{t})
\,dt+g^{k}_{t}\,dw^{k}_{t}
\end{equation}
for sufficiently large class of $f^{j},g^{k}$, say, smooth
with compact support. By the way, introducing a new
unknown function
$$
v_{t}=u_{t}-\int_{0}^{t}g^{k}_{s}\,dw^{k}_{s}
$$
reduces \eqref{6.4.1} to the heat equation with random
free term, which makes proving the {\em solvability\/}
of \eqref{6.4.1} quite elementary.
 Here is the only point
where we rely on the theory of SPDEs with
constant coefficients.

Our methods of obtaining a priori estimates
are quite different from the methods of 
\cite{Ki} and do not require developing first the
theory of SPDEs in $\bR^{d}_{+}$ or
in $\bR^{d}$ with  coefficients independent of $x$
(but {\em depending\/} on $t$ and $\omega$).
In our case
this theory  
does not help because the usual
method of freezing the coefficients
does not lead to small perturbations due to the fact
that, generally, $a$ is not continuous in~$x$.

Instead, we use new interior estimates 
of independent interest for SPDEs
in $\bR^{d}$ (Theorem \ref{theorem 12.27.2})
 which we then apply to get
an a priori estimate for equations in $\bR^{d}_{+}$
of the highest norm of the solution in terms of its lowest
norm
(Theorem \ref{theorem 12.28.1} and Corollary
\ref{corollary 12.29.1}). Then in Section
\ref{section 6.3.2} we develop a new method
of estimating the lowest norm of the solution
again avoiding considering equations
with   coefficients independent of $x$.  

We work in Sobolev spaces with weights which
is unavoidable if the stochastic terms in the
equation do not vanish  
on $\partial\bR^{d}_{+}$. It is interesting that, even if
they vanish   identically, our results are new. By the way, in that 
deterministic case
the restriction $p\geq2$ can be relaxed to $p\in(1,\infty)$
by using a standard duality argument.
Also in a standard way 
our results can be extended to cover SPDEs with VMO
coefficients in $C^{1}$ domains. The interested reader is
referred to \cite{Ki} for necessary techniques to do that.

Our results cover the classical case that $p=2$
when no continuity hypotheses is needed and
even in this case the results are new in what concerns
weights. In the case when $p\ne2$ and 
$a$ is only measurable in $x$
  the best results  
can be found in   \cite{Ki1}, where $\sigma\equiv0$ and
  $p\geq2$ is sufficiently close to 2.

\mysection{Main results}
                                       \label{section 8.27.1}

 Throughout the article
the coefficients $a^{ij}_{t}$, $a^{i}_{t}$, $b^{i}_{t}$,
$\sigma^{ik}_{t}$, $c_{t}$, and $\nu^{k}_{t}$ are 
assumed to be measurable 
with respect to $\cP\times \cB(\bR^{d})$, where $\cB(\bR^{d})$
is the Borel $\sigma$-field on $\bR^{d}$.
We understand equation 
\eqref{11.13.1}   in the sense of 
generalized functions. To be more specific
 we introduce appropriate 
Banach spaces.

Fix  some numbers
$$
p\geq2,\quad \theta\in(d-1,d-1+p), 
$$
and denote $L_{p}=L_{p}(\bR^{d}_{+})$
$$
L_{p,\theta}=\{f:M^{(\theta-d)/p}f\in L_{p}\},\quad
\|f\|_{L_{p,\theta}}=\|M^{(\theta-d)/p}f\|_{L_{p}},
$$
where $M$ is the operator of multiplying by $x^{1}$,
so that $(M^{(\theta-d)/p}f)(x)=(x^{1})^{(\theta-d)/p}f(x)$.
We use the same notation $L_{p}$ and $L_{p,\theta}$
for vector- and matrix-valued
or else
$\ell_{2}$-valued functions such as
$g_{t}=(g^{k}_{t})$ in \eqref{11.13.1}. For instance,
if $u(x)=(u^{1}(x),u^{2}(x),...)$ is 
an $\ell_{2}$-valued measurable function on $\bR^{d}$, then
$$
\|u\|^{p}_{L_{p}}=\int_{\bR^{d}_{+}}|u(x)|_{\ell_{2}}^{p}
\,dx
=\int_{\bR^{d}_{+}}\big(
\sum_{k=1}^{\infty}|u^{k}(x)|^{2}\big)^{p/2}
\,dx.
$$

 Denote
$$
D_{i}=\frac{\partial}{\partial x^{i}},\quad i=1,...,d,
\quad\Delta=D^{2}_{1}+...+D^{2}_{d}.
$$
By $Du$   we mean the gradient   with respect
to $x$ of a function $u$ on $\bR^{d}_{+}$. By $W^{1}_{p,\theta}$
we mean  the space of functions
such that $u,MDu\in L_{p,\theta}$. The norm in this space
is introduced in an obvious way. 
As is easy to see
\begin{equation}
                                                 \label{4.5.3}
\|M^{-1}u\|_{W^{1}_{p,\theta}}\sim\|M^{-1}u\|_{L_{p,\theta}}
+\|Du\|_{L_{p,\theta}}.
\end{equation}

Recall that $\tau$ is a fixed stopping time and
set
$$
\bL_{p,\theta}(\tau)=L_{p}(\opar0,\tau\cbrk,\cP,L_{p,\theta}),
\quad \bW^{1}_{p,\theta}(\tau)=
L_{p}(\opar0,\tau\cbrk,\cP,W^{1}_{p,\theta}).
$$

We also need the space $\frW^{1}_{p,\theta}(\tau)$,
which is the space of functions $u_{t}
=u_{t}(\omega,\cdot)$ on $\{(\omega,t):
0\leq t\leq\tau,t<\infty\}$ with values
in the space of generalized functions on $\bR^{d}_{+}$
and having the following properties:

(i) We have $M^{-1}u_{0}\in L_{p}(\Omega,\cF_{0},L_{p,\theta})$;

(ii)  We have $M^{-1}u
\in \bW^{1}_{p,\theta}(\tau )$;

(iii) There exist  
real valued
$f^{0}\in M^{-1}\bL_{p,\theta}(\tau)$,
$f^{1}...,f^{d}
\in  \bL_{p,\theta}(\tau)$, and an $\ell_{2}$-valued
$g=(g^{k},k=1,2,...)\in\bL_{p,\theta}(\tau)$
such that
 for any $\varphi\in C^{\infty}_{0}(\bR^{d}_{+})$ with probability 1
for all   $t\in[0,\infty)$
we have
$$
(u_{t\wedge\tau},\varphi)=(u_{0},\varphi)
+\sum_{k=1}^{\infty}\int_{0}^{t}I_{s\leq\tau}
(g^{k}_{s},\varphi)\,dw^{k}_{s}
$$
\begin{equation}
                                                 \label{1.2.01}
+
\int_{0}^{t}I_{s\leq\tau}\big(-(f^{i}_{s},D_{i}\varphi)
+(f^{0}_{s},\varphi)\big)\,ds.
\end{equation}
In particular, for any $\phi\in C^{\infty}_{0}(\bR^{d}_{+})$, 
the process
$(u_{t\wedge\tau},\phi)$ is $\cF_{t}$-adapted and continuous.

In case that property \eqref{1.2.01} holds, we write
\begin{equation}
                                       \label{12.3.1}
du_{t}=(D_{i}f^{i}_{t}+f^{0}_{t})\,dt
+g^{k}_{t}\,dw^{k}_{t}
\end{equation}
for $t\leq\tau$
and this explains the sense in which equation
\eqref{11.13.1} is understood. Of course, we still need to
specify appropriate assumptions on the coefficients
and the free terms in \eqref{11.13.1}.

For $u\in\frW^{1}_{p,\theta}(\tau)$  we write
$u\in\frW^{1}_{p,\theta,0}(\tau)$ if $u_{0}=0$. 

\begin{remark}
                                   \label{remark 4.9.1}
It is worth noting that, if $u\in\frW^{1}_{p,\theta,0}
(\tau)$, then for any $\phi
\in C^{\infty}_{0}(\bR^{d}_{+})$ the function
$u \phi\in\cW^{1}_{p,0}(\tau)$ (we remind
the definition of $\cW^{1}_{p,0}(\tau)$ later)
and as any element of $\cW^{1}_{p,0}(\tau)$ is
indistinguishable from an $L_{p}$-valued $\cF_{t}$-adapted
continuous process (see, for instance, \cite{Ito}).

\end{remark}

In the following assumption
 we use a parameter 
$K\geq0$,
which will be specified later as a small constant.

\begin{assumption} 
                                    \label{assumption 3.6.1}
 For all values of indices and arguments
we have
$$
 |Ma^{i}_{t}|+
|Mb^{i}_{t}|+| M^{2} c_{t}|
 +
| M\nu_{t} |_{\ell_{2}}\leq K ,\quad c_{t} \leq 0.
$$
\end{assumption}

\begin{remark}
Assumption \ref{assumption 3.6.1}  shows that
$a^{i}_{t},b^{i}_{t},c_{t}$, and $\nu_{t}$   go to
zero as $x^{1}\to\infty$. Actually, in applications
to SPDEs in bounded domain this is irrelevant
because far from the boundary everything is taken care
of by the theory in the whole space.
On the other hand, $a^{i}_{t},b^{i}_{t},c_{t}$, and $\nu_{t}$
  can blow  up 
to infinity for $x^{1}$ approaching zero.
\end{remark}
 
\begin{assumption}
                                    \label{assumption 3.6.2} 
For a constant  $\delta  \in(0,1]$
 for all values of the  arguments 
and    $\xi\in\bR^{d}$ we have
\begin{equation}
                                             \label{12.27.1}
a^{ij}_{t}  
  \xi^{i}
\xi^{j}\leq\delta^{-1}|\xi|^{2},\quad
(a^{ij}_{t}  
-  \alpha^{ij}_{t}) \xi^{i}
\xi^{j}\geq\delta|\xi|^{2},
\end{equation}
where
$$
\alpha^{ij}_{t}=
 (\sigma^{i\cdot} _{t},\sigma^{j\cdot}_{t})_{\ell_{2}}.
$$
\end{assumption}

Notice that we do not assume that the matrix $(a^{ij}_{t})$
is symmetric.
\begin{remark}
                                  \label{remark 4.8.1}

Observe that if $M^{-1}u\in\bW^{1}_{p,\theta}(\tau)$, then
$$
M^{-1}u\in\bL_{p,\theta}(\tau),
\quad Du\in\bL_{p,\theta}(\tau),
$$
and all
$$
 a^{ij}D_{i}u,\quad
a^{j}u,\quad Mb^{i}D_{i}u,\quad Mcu,\quad\sigma^{i}D_{i}u,
\quad\nu u
$$
belong to $\bL_{p,\theta}(\tau)$,
so that the right-hand side of \eqref{11.13.1}
has the form of the right-hand side of 
\eqref{12.3.1} with some $f^{j}$ and $g^{k}$ there
and \eqref{11.13.1} makes perfect sense for any
$u\in\frW^{1}_{p,\theta}(\tau)$.

\end{remark}
For   functions $h_{t}(x)$ on $\bR^{d+1}$
and balls $B$ in $\bR^{d}$ introduce
$$
 h_{ t(B )} 
=\frac{1}{|B |}\int_{B }h_{t}( x)\,dx,
$$
where   $|B|$ is the volume of $B$.

   If
$\rho>0$, set $B_{\rho}=\{x:|x|<\rho\}$ and for
locally integrable $h_{t}(x)$ and   continuous
 $\bR^{d}$-valued function $x_{t},t\geq0$, introduce
the integral oscillation of $h$ relative to $B$ and $x_{\cdot}$ by
$$
\text{osc}_{\rho} \, (h,x_{\cdot})=
\sup_{s\geq0}\frac{1}{\rho^{2}}
\int_{s}^{s+\rho^{2}}(|h_{t} -h_{ t(B+x_{t})}|)_{(B+x_{t})} \,dr,
$$  
where $B=B_{\rho}$.  Also for $y\in\bR^{d}$ set
$$
\text{Osc}_{\rho} \, (h,y)=
\sup_{|x_{\cdot}|_{C}\leq\rho} \sup_{r\leq\rho}\text{osc}_{r} 
\, (h,y+x_{\cdot}),
$$
where $|x_{\cdot}|_{C}$ is the sup norm of $|x_{\cdot}|$.
 Observe that $\text{osc}_{\rho} \, (h,x_{\cdot}) =0 $ if $h_{t}(x)$ is
independent of $x$.

Denote by  $\beta_{0}$   one third
of the constant $\beta_{0}(d,p,\delta)>0 $
 from Lemma 5.1 of \cite{Kr09}.
\begin{assumption}
                                  \label{assumption 4.7.1}

There exist a constant  $\varepsilon
\in(0,1]$ 
such that for any   $y\in\bR^{d}_{+}$
(and $\omega$) 
we have 
\begin{equation}
                                                   \label{8.7.1}
\text{Osc}_{\varepsilon y^{1}} \,(a^{ij},y)\leq \beta_{0},
\quad\forall i,j.
\end{equation}

Furthermore,
$$
(a^{jk}_{t}(x)  
-  \alpha^{jk}_{t}(y)) \xi^{j}\xi^{k}\geq\delta|\xi|^{2}
$$
for all   $t $, $\xi$, and $x$  satisfying
  $|x-y|\leq\varepsilon y^{1}  $.

\end{assumption}   

\begin{remark}
                                            \label{remark 6.3.1}

This assumption is quite substantially weaker than
similar assumptions known in the literature
(see, for instance, \cite{Ki} and the references therein),
where the oscillation of $a^{ij}$ in \eqref{8.7.1} 
is understood as
\begin{equation}
                                                     \label{8.27.1}
\sup_{t\geq0}\sup_{|x-y|\leq\varepsilon(x^{1}\wedge
y^{1})}|a^{ij}_{t}(x)-a^{ij}_{t}(y)|.
\end{equation}
 
 It is easy to see that if, for an  $\varepsilon\in(0,1]$,
 \eqref{8.27.1} is less than a $\beta>0$,
then the left hand-side of \eqref{8.7.1} is also less than $\beta$
if we replace there $\varepsilon$ with $\varepsilon/4$.) 
With such substitution $a^{ij}_{t}(x)$ 
will have jumps at each point $x\in\bR^{d}_{+}$ not larger than
$\beta_{0}$, which is a small constant.

On the other hand,
if $a^{ij}_{t}(x)$
is independent of $t$,  
then, for $0< y^{1} \leq2$,  \eqref{8.7.1} 
 is satisfied if $a\in \text{VMO}$, which
is the class of functions with vanishing mean oscillation
and which for $d=2$ contains, for instance,
the function $2+\sin f(x)$, where $f(x)=\ln^{1/3}(|x-e|\wedge1)$
and $e$ is the first basis vector in $\bR^{d}$. The 
usual oscillation of this function at $e$ is 2.

\end{remark}
 \begin{remark}
It follows from our proofs that if $\sigma\equiv0$, then
we can relax condition \eqref{8.7.1} by using the modified 
  integral oscillations which are defined by taking
$x_{t}\equiv0$.
\end{remark}
Let $\beta_{1}=\beta_{1}(d,p,\delta,\varepsilon )>0$
be the constant   
from  Lemma 5.2 of \cite{Kr09}.

\begin{assumption} 
                                       \label{assumption 8.10.1}
There exists a constant  $\varepsilon_{1}>0$ 
such that for any $t\geq0$  
we have
$$
|\sigma^{i\cdot}_{t}(x)-\sigma^{i\cdot}_{t}(y)
|_{\ell_{2}}\leq \beta_{1},
$$
whenever 
$$
x,y\in\bR^{d}_{+}, 
\quad |x-y|\leq 
\varepsilon_{1}(x^{1}\wedge y^{1}),\quad i=1,...,d . 
$$

\end{assumption}

Our first main result is the following.

\begin{theorem}
                                   \label{theorem 5.30.1}
Let $\bar{\delta}>0$ be  a constant such that
for any $\xi\in\bR^{d}$ and all values of arguments
we have
\begin{equation}
                                                  \label{6.1.1}
\bar{\delta}\big(\sum_{i}
a^{i1}_{t}\xi^{i}\big)^{2}\leq  (a^{ij}_{t}-
\alpha^{ij}_{t})\xi^{i}\xi^{j}.
\end{equation}
Let  Assumptions \ref{assumption 3.6.1} through 
 \ref{assumption 8.10.1}
  be satisfied 
 with a (small) constant $K=K(d,p,\delta,
\bar{\delta},\theta,\varepsilon,\varepsilon_{1})>0$,
an estimate from below for which can be obtained
from the proof.        
Set
$$
\gamma=\theta-d-p+1\quad(<0)
$$
and assume that
$$
|\gamma|(\bar{\delta}\delta)^{1/2}(p-1)>p|\gamma+1|,
$$
which holds, for instance, if 
$
\theta=d+p-2 $ when $\gamma+1=0$.
Then   for any  $f^{0}, ..., f^{d}$, and
$g=(g^{k})$ satisfying
$$
Mf^{0}, f^{i},
g=(g^{k})\in\bL_{p,\theta}(\tau),\quad i=1,...,d
$$
there exists a unique
$u
\in\frW^{1}_{p,\theta,0}(\tau)$ satisfying \eqref{11.13.1}
in $\bR^{d}_{+}$.
Furthermore, for this solution
\begin{equation}
                                               \label{5.30.1}
\| M^{-1}u\|_{\bW^{1}_{p,\theta}(\tau)}  
\leq N\big(\| Mf^{0}\|_{\bL_{p,\theta}(\tau)}
+\sum_{i=1}^{d}\|  f^{i}\|_{\bL_{p,\theta}(\tau)}
+ \| g\|_{\bL_{p,\theta}(\tau)}\big),
\end{equation}
where $N$ depends only on
$d,p,\delta,\theta,\bar{\delta}$, $\varepsilon $, and 
$\varepsilon_{1}$.
\end{theorem}

\begin{remark}
                                           \label{remark 6.4.1}
As it follows from
the proof of Theorem \ref{theorem 5.30.1}, if $p=2$,
Assumptions \ref{assumption 4.7.1} and 
\ref{assumption 8.10.1}
are not needed.
Thus we obtain the classical result
on Hilbert space solvability of SPDEs in half spaces
with one improvement that we can allow spaces with weights.
By the way, observe that the proof of Theorem \ref{theorem 5.30.1}
does not use  the Hilbert space theory of SPDEs.

\end{remark}

To state our second result we need an additional
assumption.

\begin{assumption}
                                     \label{assumption 6.2.1}

(i) There exists a constant $\tilde{\delta}
\in(0,1]$ such that for all   $\xi\in\bR^{d}$
and all arguments
we have
\begin{equation}
                                             \label{2.23.6}
\tilde{\delta} (\sum_{j}\hat{a}^{1j}_{t}
\xi^{j})^{2}\leq a^{11}_{t}(a^{ij}_{t}-\alpha^{ij}_{t}) \xi^{i}\xi^{j},
\end{equation}
where
$$
\hat{a}^{ij}_{t}=(1/2)(a^{ij}_{t}+a^{ji}_{t}).
$$

(ii) It holds that
\begin{equation}
                                       \label{1.24.5}
d-1+p\big[1-\frac{1}{p(1-\tilde{\delta})
+\tilde{\delta}}\big]<
\theta<d-1+p.
\end{equation}

(iii) For a constant  $ \beta_{2}>0$, if
 $x,y\in\bR^{d}_{+}$ are such that 
$|x-y|\leq 
  x^{1}\wedge y^{1} $, then for all
$ i =1,...,d$ and  $t>0$
\begin{equation}
                                       \label{6.3.1}
|\hat{a}^{i1}_{t}(x)-\hat{a}^{i1}_{t}(y)|\leq \beta_{2}.
\end{equation}

\end{assumption}

\begin{remark}
In previous works on a similar subject (see, for instance,
\cite{Ki} or \cite{KL}) 
a  condition stronger than \eqref{2.23.6} used to
be assumed:
\begin{equation}
                                       \label{6.2.5}
\tilde{\delta}a^{ij}_{t}\xi^{i}\xi^{j}\leq 
(a^{ij}_{t}-\alpha^{ij}_{t}) \xi^{i}\xi^{j}.
\end{equation}
That \eqref{6.2.5} is stronger than \eqref{2.23.6}
follows from the fact that  
 for the positive definite matrix $(\hat{a}^{ij}_{t})$
and $\eta=(1,0,...,0)$
it holds that
$$
 (\sum_{j}\hat{a}^{1j}_{t}
\xi^{j})^{2}=
 (\sum_{j}\hat{a}^{ij}_{t}
\eta^{i}\xi^{j})^{2}\leq(\hat{a}^{ij}_{t}\eta^{i}\eta^{j})
\hat{a}^{ij}_{t}\xi^{i}\xi^{j}=a^{11}_{t}
a^{ij}_{t}\xi^{i}\xi^{j}.
$$
Also observe that
sometimes \eqref{2.23.6} holds
with $\tilde{\delta}=1$ and \eqref{6.2.5} does not.
 This happens, for instance,
if $\alpha^{1j}_{t}\equiv0$
for all $j$ and $\hat{a}^{1j}_{t}\equiv0$
for $j\ne1$. Finally, in the   case
when $\sigma\equiv0$ condition \eqref{2.23.6}
is satisfied with $\tilde{\delta}=1$ and then condition
\eqref{1.24.5} becomes $d-1<\theta<d-1+p$ which is the widest
range possible for $\theta$ even in the deterministic case
for the heat equation.
\end{remark}
\begin{remark}
Condition \eqref{6.3.1} is imposed only on
$\hat{a}^{i1}_{t}$. As is discussed in \cite{Ki}
(also see references therein), this condition
allows rather sharp oscillations of $ \hat{a}^{i1}_{t}(x) $
near $\partial\bR^{d}_{+}$. The other 
entries of $(a^{ij}_{t}(x))$
are still allowed to be discontinuous in $x$ but yet kind of
belong to VMO (cf. Remark \ref{remark 6.3.1}).

\end{remark}

\begin{theorem}
                                   \label{theorem 6.2.1}

  There exist (small) constants $K>0$ and $\beta_{2}>0$,
depending only on
$d,p,\delta,
\tilde{\delta},\theta,\varepsilon$, and $\varepsilon_{1} $ and
  estimates from below for which can be obtained
from the proof, such that if
Assumptions \ref{assumption 3.6.1}  
through \ref{assumption 6.2.1}  are satisfied
with these constants, then
 the assertion of
Theorem \ref{theorem 5.30.1} holds true again
with $\tilde{\delta} $ in place of
$\bar{\delta}$ in the arguments of $N$. 
\end{theorem}

We prove Theorems \ref{theorem 5.30.1} and \ref{theorem 6.2.1}
in Section \ref{section 6.4.2} after preparing
necessary tools in Section  \ref{section 6.4.3},
where we treat equations in $\bR^{d}$,
and in Section \ref{section 6.3.2} containing
auxiliary results for equations in $\bR^{d}_{+}$.

\mysection{Auxiliary results for equations in $\bR^{d}$}
                                       \label{section 6.4.3}

The assumptions in this section are somewhat different
from the assumptions of Section \ref{section 8.27.1}
apart from the assumption about the measurability
of the coefficients. 

To investigate the equations in $\bR^{d}_{+}$
we need a few results about equations in $\bR^{d}$.
  To state them
we remind the reader  the definition
of spaces $\bW^{1}_{p}(\tau)$ and $\cW^{1}_{p}(\tau)$
introduced in \cite{Ito} (which is somewhat different
from $\cH^{1}_{p}(\tau)$
in \cite{Ki} or \cite{Kr99}, see the discussion of the differences
in \cite{Kr09}).

As usual,  
$$
W^{1}_{p}=\{u\in L_{p}(\bR^{d}): Du\in L_{p }(\bR^{d})\},
\quad
 \|u\|_{W^{1}_{p}}=
\|u\|_{L_{p}(\bR^{d})}+\|Du\|_{L_{p}(\bR^{d})}.
$$

Recall that $\tau$ is a stopping time and set
$$
\bL _{p}(\tau):=L_{p}(\opar 0,\tau\cbrk,\cP,
L_{p}(\bR^{d})),\quad
\bW^{1}_{p}(\tau):=L_{p}(\opar 0,\tau\cbrk,\cP,
W^{1}_{p}).
$$
The space $\cW^{1}_{p}(\tau)$,
  is introduced as the space of functions $u_{t}
=u_{t}(\omega,\cdot)$ on $\{(\omega,t):
0\leq t\leq\tau,t<\infty\}$ with values
in the space of generalized functions on $\bR^{d}$
and having the following properties:

(i) We have $u_{0}\in L_{p}(\Omega,\cF_{0},L_{p}(\bR^{d}))$;

(ii)  We have $u
\in \bW^{1}_{p}(\tau )$;

(iii) There exist   $f^{i}\in \bL_{p}(\tau)$,
$i=0,...,d$, and $g=(g^{1},g^{2},...)\in \bL_{p}(\tau)$
such that
 for any $\varphi\in C^{\infty}_{0}(\bR^{d})$ 
with probability 1
for all   $t\in[0,\infty)$
we have
$$
(u_{t\wedge\tau},\varphi)=(u_{0},\varphi)
+\sum_{k=1}^{\infty}\int_{0}^{t}I_{s\leq\tau}
(g^{k}_{s},\varphi)\,dw^{k}_{s}
$$
\begin{equation}
                                                 \label{1.2.1}
+
\int_{0}^{t}I_{s\leq\tau}\big((f^{0}_{s},\varphi)-(f^{i}_{s},D_{i}\varphi)
 \big)\,ds.
\end{equation}
In particular, for any $\phi\in C^{\infty}_{0}$, the process
$(u_{t\wedge\tau},\phi)$ is $\cF_{t}$-adapted and continuous.

The following result is 
a somewhat weakened version of Corollary 5.5 in~\cite{Kr09}.

\begin{lemma}
                                    \label{lemma 12.7.1}
Let $G\subset\bR^{d}$ be a domain (perhaps, $G=\bR^{d}$)
and let $K\geq0$, $\varepsilon >0$, and $\varepsilon_{1}
\in(0,\varepsilon/4]$ be some constants.

(i) Let  
$f^{j},g\in\bL_{p}(\tau)$, and
 $u\in\cW^{1}_{p,0}(\tau)$
satisfy~\eqref{11.13.1} in $\bR^{d}$  for $t\leq\tau$
and be such that
$u_{t}(x)=0$ if $x\not\in G$.

(ii) Suppose that Assumption
\ref{assumption 3.6.2}
is  satisfied and suppose that for $y\in G$
and all values of indices and other arguments
$$
 |a^{i}_{t}(y)|+|b^{i}_{t}(y)|+|c_{t}(y)|
+|\nu_{t}(y)|_{\ell_{2}}\leq
K,\quad c_{t}(y)\leq0.
$$

(iii) Assume that,
for any  $x_{0} $, such 
that $\text{\rm dist}\,(x_{0},G)\leq\varepsilon_{1}$,
we have
\begin{equation}
                                             \label{8.12.1}
\text{\rm Qsc}_{\varepsilon}\,(a^{ij},x_{0})\leq
\beta_{0},\quad\forall i,j,
\end{equation}
where $\beta_{0}$ is the one third
of $\beta_{0}(d,p,\delta)>0 $
 from Lemma 5.1 of \cite{Kr09}, and
\begin{equation}
                                             \label{8.12.2}
|\sigma^{i\cdot}_{t}(x)-
\sigma^{i\cdot}_{t}(x_{0})|_{\ell_{2}}\leq\beta_{1},\quad 
(a^{jk}_{t}(y)  
-  \alpha^{jk}_{t}(x_{0})) \xi^{j}
\xi^{k}\geq\delta|\xi|^{2}
\end{equation}
for all values of indices and arguments 
such that $|x-x_{0}|\leq\varepsilon_{1}$ and
$|y-x_{0}|\leq\varepsilon$, where $\beta_{1}=
\beta_{1}(d,\delta,p,\varepsilon/2)>0$  is the constant
from  Lemma 5.2 of \cite{Kr09}.

Then there exist a constant   $N $ 
depending only on $d,p,K,\delta $, $\varepsilon$, 
  and  $\varepsilon_{1}$ such that  
 $$
 \|Du \|_{\bL_{p}(\tau)} 
\leq N\big(\sum_{i=1}^{d}
\|f^{i}\|_{\bL_{p}(\tau)}+\|g\|_{\bL_{p}(\tau)}
+\|f^{0} \|^{1/2}_{\bL_{p}(\tau)}\|u\|_{\bL_{p}(\tau)}^{1/2}
+\|u\|_{\bL_{p}(\tau)}
\big).
$$
\end{lemma}

Next we give a version of Lemma \ref{lemma 12.7.1} 
for some particular domains $G$ the most important of
which will be $\{|x^{1}|\leq R\}$.
We state it in a slightly more general setting
suitable for investigating interior smoothness
of solutions in $\bR^{d}$ or in $\bR^{d}_{+}$.

We fix an integer $d_{1}\in[1,d]$ and for $x\in\bR^{d}$
introduce 
$$
|x|'=\big(\sum_{i=1}^{d_{1}}(x^{i})^{2}
\big)^{1/2} ,\quad B'_{R}=\{x\in\bR^{d}:
|x|'<R\}.
$$

\begin{theorem}  
                                         \label{theorem 12.29.1}
Take some $\varepsilon>0$, $\varepsilon_{1}\in(0,\varepsilon/4]$,
$K\geq0$, and $R>0$.

(i) Let  
$f^{j},g\in\bL_{p}(\tau)$, and
 $u\in\cW^{1}_{p,0}(\tau)$
satisfy  \eqref{11.13.1} in $\bR^{d}$ for $t\leq\tau$.

(ii) Assume that  $u_{t}(x)=0$ if $x\not\in B'_{R}$.

(iii) Suppose that Assumption 
\ref{assumption 3.6.2}
is  satisfied and for $y\in B'_{R}$ and
all values of the indices and
other arguments
$$
R|a^{i}_{t}(y)|+ R|b^{i}_{t}(y)|+
R|\nu_{t}(y)|_{\ell_{2}}+R^{2}|c_{t}(y)| \leq K,
\quad c_{t}(y)\leq0.
$$

(iv)   Assume that \eqref{8.12.1} 
with $\varepsilon R$ in place of 
$\varepsilon$ and \eqref{8.12.2}
hold for any
$x_{0}$, such that $|x _{0}|'\leq (1+\varepsilon)R$,
and $x,y$ such that
$|x-x_{0}|\leq\varepsilon_{1}R$, $|y-x_{0}|\leq
\varepsilon R,
$
and all values of indices and other arguments.

Then there exists a constant 
  $N =N (d,p,\delta ,K,\varepsilon,\varepsilon_{1})$
such that
\begin{equation}
                                                \label{12.29.2}
\| Du\|_{\bL_{p}(\tau)}
\leq N\big(\sum_{i=1}^{d}\|f^{i}\|_{\bL_{p}(\tau)}
+  \| g\|_{\bL_{p}(\tau)} 
+   \|f^{0}\|_{\bL_{p}(\tau)}^{1/2}
\|u\|_{\bL_{p}(\tau)}^{1/2}+ R^{-1} \|u\|_{\bL_{p}(\tau)}\big).
\end{equation}

\end{theorem}

Proof. If $R=1$, the result follows directly from 
Lemma~\ref{lemma 12.7.1}.
The case of general $R$ we reduce to the particular one
by using dilations. Introduce
$$
\hat\cF _{t}=\cF_{R^{2}t},\quad\hat\tau =R^{-2}\tau,\quad
\hat w^{ k}_{t}=R^{-1}w^{k}_{R^{2}t},
$$
$$
(\hat a^{ij} _{t},\hat{a}^{i}_{t} ,\hat b _{t} ,\hat c _{t} ,
\hat\sigma _{t},\hat\nu _{t})(x)
=(a^{ij} _{R^{2}t},Ra^{i}_{R^{2}t} ,
Rb _{R^{2}t} ,R^{2}c _{R^{2}t} ,
\sigma  _{R^{2}t},R\nu _{R^{2}t})(Rx),
$$
$$
\hat u _{t}(x)=u_{R^{2}t}(Rx),\quad
\hat f^{i}_{t}(x)=Rf^{i}_{R^{2}t}
(Rx),\quad i=1,...,d,
$$
$$
\hat f^{0}_{t}(x)=R^{2}f^{0}_{R^{2}t},\quad 
\hat g^{k}_{t}(x)=Rg^{k}_{R^{2}t}(Rx).
$$
Also introduce the operators $\hat L _{t}$ and
$\hat\Lambda^{ k}_{t}$ constructing them from the above introduced
coefficients.
It is easily seen
 that $\hat w^{k}_{t}$ are independent $\hat\cF_{t}$-Wiener
processes, $\hat\tau$ is an $\hat\cF_{t}$-stopping time,
  all the above processes with hats
are predictable with respect to  the
filtration $\{\hat\cF_{t}\}$, and
$\hat{u}\in\hat\cW^{1}_{p,0}(\hat\tau)$,
$\hat f,\hat g\in\hat{\, \bL}_{p}(\hat\tau)$, 
 where the spaces with hats
are defined on the basis of $\{\hat\cF_{t}\}$.

Observe that for $t<\hat\tau$
$$
\hat L _{t}\hat u _{t}(x)=\big(D_{j}(
\hat a^{ ij}(x)D_{i }\hat u _{t}(x)+\hat{a}^{j}_{t}
(x)\hat{u}_{t}(x))  
+\hat b^{ i}_{t}(x)D_{i}\hat u_{t}(x)+\hat c_{t}(x)\hat u_{t}(x)
$$
$$
=R^{2} \big(D_{j}( a^{ij}_{R^{2}t}D_{i }u_{R^{2}t}
+a^{j}_{R^{2}t}u_{R^{2}t})+b^{i}_{R^{2}t}
D_{i}u_{R^{2}t}+c_{R^{2}t}u_{R^{2}t}\big)(Rx)
$$
$$
=R^{2}L_{R^{2}t}u_{R^{2}t}(Rx),\quad D_{i}\hat{f}^{i}_{R^{2}t}
(x)=R^{2}(D_{i}f ^{i}_{R^{2}t})(Rx),
$$
$$
\int_{0}^{t\wedge\hat{\tau}}[\hat L _{s}\hat 
u _{s}(x)+D_{i}\hat f^{i}_{s}(x)
+\hat{f}^{0}_{s}(x)]\,ds
$$
$$
=\int_{0}^{(R^{2}t)\wedge \tau }[ L _{s}u _{s}(Rx)+ D_{i}f^{i}_{s}(Rx)
+f^{0}_{s}(Rx)]\,ds.
$$
Of course, we understand this equality in the sense
of distributions:
$$
\int_{0}^{t\wedge\hat{\tau}}(\hat L _{s}\hat u _{s} +
D_{i}\hat f^{i}_{s}+\hat{f}^{0}_{s} ,\phi)\,ds
=\int_{0}^{(R^{2}t)\wedge \tau }([ L _{s}u _{s}+ 
D_{i}f^{i}_{s}+f^{0}_{s}](R\cdot),\phi)\,ds
$$
for any $\phi\in C^{\infty}_{0}(\bR^{d})$.
 One also knows that
if $\hat h_{t}$ is an $\hat\cF_{t}$-predictable process
satisfying a natural integrability condition with respect to $t$,
then
$$
\int_{0}^{t}\hat h_{s}\,d\hat w^{k}_{s}
=R^{-1}\int_{0}^{R^{2}t}\hat h_{R^{-2}s}\,dw^{k}_{s} 
\quad\text{(a.s.)}.
$$
Therefore, (a.s.)  
$$
\int_{0}^{t\wedge\hat{\tau}}[\hat\Lambda^{k}_{s}\hat u_{s}+
 \hat g^{k}_{s}](x)\,d\hat w^{k}_{s}
=R\int_{0}^{t\wedge\hat{\tau}}[ \Lambda^{k}_{R^{2}s}  u_{R^{2}s}+
   g^{k}_{R^{2}s}](Rx)\,d\hat w^{k}_{s}
$$
$$
=\int_{0}^{(R^{2}t)\wedge\tau}[ \Lambda^{k}_{s}  u_{ s}+
   g^{k}_{ s}](Rx)\,d  w^{k}_{s}.
$$

It follows that (a.s.)  
$$
\int_{0}^{t\wedge\hat{\tau}}[\hat L _{s}\hat u _{s}(x)+
D_{i}\hat f^{i}_{s}(x)
+\hat{f}^{0}_{s}(x)]\,ds
$$
$$
+\int_{0}^{t\wedge\hat{\tau}}[\hat\Lambda^{k}_{s}\hat u_{s}+
 \hat g^{k}_{s}](x)\,d\hat w^{k}_{s}=u_{(R^{2}t)\wedge\tau}(Rx)=
\hat u_{t\wedge\hat{\tau}}(x),
$$
so that $\hat u$ satisfies equation \eqref{11.13.1}
with new operators and free terms. It is also easy to see that
our objects with hats satisfy the assumptions of the theorem
with $R=1$. Therefore, by the result for $R=1$
$$
 \|D \hat u\|_{\hat{\,\bL}_{p}(\hat\tau)}
\leq N\big(\sum_{j=1}^{d}
\|\hat f^{j}\|_{\hat{\,\bL}_{p}(\hat\tau)}+
\|\hat g\|_{\hat{\,\bL}_{p}(\hat\tau)}+
\|\hat f^{0}\|_{\hat{\,\bL}_{p}(\hat\tau)}^{1/2}
\|\hat u\|_{\hat{\,\bL}_{p}(\hat\tau)}^{1/2}+
\|\hat u\|_{\hat{\,\bL}_{p}(\hat\tau)} \big).
$$
Now it only remains to notice that changing variables shows that
this inequality is precisely \eqref{12.29.2}.
The theorem is proved.

Here is an interior estimate for equations in $\bR^{d}$.
In its spirit it is similar to Theorem 2.3 of \cite{Kr08}.

\begin{theorem}
                                \label{theorem 12.27.2}
Let
assumptions  (i), (iii), and (iv)  of Theorem \ref{theorem 12.29.1}
be satisfied. 
Then, for any  $r\in(0,R)$, we have 
$$
\|I_{B'_{r}}D u\|_{\bL_{p}(\tau)}\leq N\big(
 \|I_{B'_{R}}f^{0}\|_{\bL_{p}(\tau)}^{1/2}
\|I_{B'_{R}}u\|_{\bL_{p}(\tau)}^{1/2}
$$
\begin{equation}
                                                 \label{1.1.2}
+
\sum_{i=1}^{d}\|I_{B'_{R}}f^{i}\|_{\bL_{p}(\tau)}+
\|I_{B'_{R}}g\|_{\bL_{p}(\tau)}\big)
+N(R-r)^{-1}\|uI_{B'_{R}}\|_{\bL_{p}(\tau)},
\end{equation}
where $N=N( d,p,\delta,K,\varepsilon ,\varepsilon_{1})$.

\end{theorem}

Proof. We follow a usual procedure taken
from the theory of PDEs.
Let $\chi(s)$ be an infinitely differentiable function
on $\bR$ such that $\chi(s)=1$ for $s\leq0$ and $\chi(s)=0$
for $s\geq1$. For $m=0,1,2,...$ introduce, ($r_{0}=r$)
$$
r_{m}=r+(R-r)\sum_{j=1}^{m}2^{-j},\quad
\zeta_{m}(x)=\chi\big(2^{m+1}(R-r)^{-1}(|x|'-r_{m}) \big).
$$
As is easy to check, for
$$
Q(m)=  B'_{r_{m}}\,,
$$
it holds that
$$
\zeta_{m}=1\quad\text{on}\quad Q(m),\quad\zeta_{m}=0\quad
\text{outside}
\quad  Q(m+1).
$$
Also (observe that $N2^{m+1}=N_{1}2^{m}$ with $N_{1}=2N$) 
$$
|D\zeta_{m}|\leq N2^{m}(R-r)^{-1} .
$$

Next, the function $\zeta_{m}u_{t}$ is in
$\cW^{1}_{p,0}(\tau)$ and satisfies 
\begin{equation}
                                           \label{12.11.3}
d(\zeta_{m}u _{t})=\big(L_{t}(\zeta_{m}u _{t}) 
+D_{j} f ^{ j}_{mt}+ f ^{ 0}_{mt}\big)\,dt
+\big(\Lambda_{t}^{k}(\zeta_{m}u_{t})
+ g ^{k }_{mt}\big)\,dw^{k}_{t},
\end{equation}
where  
$$
 f ^{ j}_{mt}= -
a^{ij}_{t} u_{t}D_{i}\zeta_{m}
  +\zeta_{m}f^{j}_{t},\quad j=1,...,d,
$$
$$
 f ^{ 0}_{mt}=-a^{ij}_{t}(D_{i}u_{t})D_{j}\zeta_{m}
-u_{t}a^{j}_{t}D_{j}\zeta_{m}
- u_{t}b^{i}_{t}D_{i}\zeta_{m} +\zeta_{m}f^{0}_{t} -
f^{i}_{t}D_{i}\zeta_{m}.
$$
$$
 g ^{k}_{mt}= 
\zeta_{m}g^{k}_{t} 
-u_{t}\sigma^{ik}_{t}D_{i}\zeta_{m}.
$$

Notice that
$$
| f ^{ j}_{mt}|\leq N 2^{m}(R-r)^{-1}|\zeta_{m+1}  u_{t}|
+ \zeta_{m}|f^{j}_{t}|,\quad j=1,...,d,
$$
$$
| f ^{ 0}_{mt}|\leq N 2^{m}(R-r)^{-1}\zeta_{m+1}|Du_{t}|
+N 2^{m}(R-r)^{-2}|\zeta_{m+1} u_{t}| 
$$
$$
+N \zeta_{m}|f^{0}_{t}|+
N2^{m}(R-r)^{-1}\zeta_{m+1}\sum_{j=1}^{d}|f^{j}_{t}|
$$
$$
\leq N 2^{m}(R-r)^{-1}|D(\zeta_{m+1}u_{t})|
+N 2^{2m}(R-r)^{-2}|\zeta_{m+1} u_{t}| 
$$
$$
+N \zeta_{m}|f^{0}_{t}|+
N2^{m}(R-r)^{-1}\zeta_{m+1}\sum_{j=1}^{d}|f^{j}_{t}|,
$$
$$
|g _{mt}|_{\ell_{2}}\leq 
 \zeta_{m}|g_{t}|_{\ell_{2}}
+N 2^{m}(R-r)^{-1}|\zeta_{m+1}u_{t}|.
$$
 
Since $\zeta_{m}u_{t}(x)=0$ for $x\not\in B'_{R}$,
by Theorem \ref{theorem 12.29.1} and Young's inequality 
  we have
$$
D_{m}:= \|D (\zeta_{m}u )
\|_{\bL_{p}(\tau)} 
  \leq NF
+N 2^{m}(R-r)^{-1} U_{m+1}
$$
$$
+N2^{m/2}(R-r)^{-1/2}D_{m+1}^{1/2}U^{1/2}_{m+1}
$$
$$
  \leq NF
+N 2^{m}(R-r)^{-1} U_{m+1}+2^{-2}D_{m+1},
$$
where
$$
U_{m }:=\| \zeta_{m}u 
\|_{\bL_{p}(\tau)} ,\quad
F:=
\sum_{i=1}^{d}\|I_{B'_{R}}f^{i}\|_{\bL_{p}(\tau)}+
\|I_{B'_{R}}g\|_{\bL_{p}(\tau)}
$$
$$
+\|I_{B'_{R}}f^{0}\|_{\bL_{p}(\tau)}^{1/2}
\|I_{B'_{R}}u\|_{\bL_{p}(\tau)}^{1/2}.
$$
It follows that
$$
D_{0}+\sum_{m=1}^{\infty}2^{-2m}D_{m}\leq
N F+N(R-r)^{-1}\|uI_{B'_{R}}\|_{\bL_{p}(\tau)}
+\sum_{m=1}^{\infty}2^{-2m}D_{m}.
$$
By canceling like terms we estimate $D_{0}$ by the right-hand
side of \eqref{1.1.2}. Its left-hand side is certainly smaller than
$D_{0}$. This would yield \eqref{1.1.2} provided that
what we canceled is finite.

Obviously, 
$$
D_{m}\leq N\|Du\|_{\bL_{p}(\tau)}+
N2^{m}(R-r)^{-1}\|u\|_{\bL_{p}(\tau)}
$$
and the terms in question are finite since $u\in\bW^{1}_{p}
(\tau)$.   The theorem is proved.

\mysection{Auxiliary results for 
equations in $\bR^{d}_{+}$}
                                            \label{section 6.3.2}

In this section we are  investigating the local regularity
of solutions in $\bR^{d}_{+}$ and give preliminary
a priori estimates.

For $r>0$ denote
$$
G_{r}=\{x\in\bR^{d}:0<x^{1}<r\}.
$$
Here is the divergence form counterpart of Theorem 4.3
of \cite{Kr08}.
 \begin{theorem}
                                      \label{theorem 12.28.1}

Take an $R\in(0,\infty]$ and  
  suppose the following.

(i)  Assumptions \ref{assumption 3.6.1} through 
 \ref{assumption 8.10.1}
  are satisfied.

(ii) We have a function $u$ such that 
$\phi u\in \cW^{1}_{p,0}(\tau)$ for any
$\phi\in C^{\infty}_{0}(G_{R})$ and $ u$
satisfies \eqref{11.13.1} in $\bR^{d}_{+}$ for $t\leq\tau$
with some   
$
f^{j},g=(g^{k},k=1,2,...)
$
 such that
$Mf^{0},f^{i},g\in\bL
_{p,\theta}(\tau)$, $i=1,...,d$.

Then, for any $r\in(0,R/4)$,
$$
\|I_{G_{r}} D u\|_{\bL_{p,\theta}(\tau)}
\leq
N\|I_{G_{R}}Mf^{0}\|_{\bL_{p,\theta}(\tau)}^{1/2}
\|I_{G_{R}}M^{-1}u\|_{\bL_{p,\theta}(\tau)}^{1/2}
$$
\begin{equation}
                                               \label{1.1.4}
+ N\sum_{i=1}^{d}\|I_{G_{R}}f^{i}\|_{\bL_{p,\theta}(\tau)}
+N\|I_{G_{R}}g\|_{\bL_{p,\theta}(\tau)}
+N\|I_{G_{R}}M^{-1}u\|_{\bL_{p,\theta}(\tau)},
\end{equation}
where   $N=N(d,p,\delta ,\varepsilon, \varepsilon_{1},K)$.
\end{theorem}
 
Proof. We are going to apply Theorem \ref{theorem 12.27.2}
to shifted $B'_{R}$ when $d_{1}=1$.
For $n=-1,0,1,...$, set $r_{n}=2^{-n/3}r$. Observe that
if $n\geq0$, then the half width of
$G_{r_{n-1}}\setminus G_{r_{n+2}}$ equals $\rho_{n}:= r_{n+2}/2$
 and 
$$
r_{n-1}+ \rho_{n}\leq 2r_{-1} <4r<R,
\quad
r_{n+2}- \rho_{n}=\rho_{n} .
$$
Let $c_{n}=(r_{n-1}+r_{n+2})/2$ and observe that
for $x_{0}\in\bR^{d}_{+}$ such that $|x_{0}^{1}
-c_{n}|\leq(1+\varepsilon)\rho_{n}$ we have
$\rho_{n}\leq x_{0}^{1}$ because $\varepsilon\leq1 $.
It follows that
$$
\text{Osc\,}_{\varepsilon\rho_{n}}(a^{ij},x_{0})\leq\beta_{0}.
$$
Also, for $y\in G_{r_{n-1}}\setminus G_{r_{n+2}}$
we have $\rho_{n}\leq y^{1}$ and
$$
\rho_{n}|a^{i}_{t}(y)|+ \rho_{n}|b^{i}_{t}(y)|+
\rho_{n}|\nu_{t}(y)|_{\ell_{2}}+\rho_{n}^{2}|c_{t}(y)| \leq K,
\quad c_{t}(y)\leq0.
$$

Next, if $|x_{0}^{1}
-c_{n}|\leq(1+\varepsilon)\rho_{n}$  
  and $|y-x_{0}|\leq
\varepsilon\rho_{n}$, then
$|y-x_{0}|\leq
\varepsilon x_{0}^{1}$ and
$$
(a^{jk}_{t}(y)  
-  \alpha^{jk}_{t}(x_{0})) \xi^{j}
\xi^{k}\geq\delta|\xi|^{2}.
$$
Finally, define $\gamma\in(0,\varepsilon/4]$ by
$$
\frac{\gamma}{1-\gamma}=\varepsilon_{1}\wedge\frac{\varepsilon}{4}
$$
and observe that
if $|x_{0}^{1}
-c_{n}|\leq(1+\varepsilon)\rho_{n}$ and $|x-x_{0}|\leq
\gamma\rho_{n}$, then
\begin{equation}
                                           \label{8.27.3}
|x-x_{0}|\leq
\gamma x_{0}^{1}\leq\gamma (x_{0}^{1}\wedge x^{1})
\leq\varepsilon_{1}(x_{0}^{1}\wedge x^{1})
\end{equation}
if $x^{1}\geq x_{0}^{1}$ and,  if $x^{1}< x_{0}^{1}$,
then $x_{0}^{1}-x^{1}\leq\gamma x_{0}^{1}$,
$x_{0}^{1}\leq(1-\gamma)^{-1}(x_{0}^{1}\wedge x^{1})$
and the inequality between the extreme terms in \eqref{8.27.3}
holds again. In that case 
$$
|\sigma^{i\cdot}_{t}(x)-
\sigma^{i\cdot}_{t}(x_{0})|_{\ell_{2}}\leq\beta_{1}.
$$
This means that the assumptions of Theorem \ref{theorem 12.27.2}
about the coefficients are satisfied if we shift
$c_{n}$ into the origin.

Furthermore, if $n\geq0$, $\zeta\in C^{\infty}_{0}((0, R))$, and
$\zeta(z)=1$ for $r_{n+2}\leq z\leq r_{n-1}$, then 
$\zeta u$ satisfies  \eqref{11.13.1} in $\bR^{d}  $
with certain $f$ and $g$ which on $  G_{r_{n-1}}
\setminus G_{r_{n+2}}$
coincide with the original ones. Finally,
if $n\geq0$, then the distance between the boundaries of
$G_{r_{n}}\setminus G_{r_{n+1}}$ and 
$G_{r_{n-1}}\setminus G_{r_{n+2}}$ is $(2^{1/3}-1)r_{n+2}$. 
 
It follows by Theorem
\ref{theorem 12.27.2} that for $n\geq0$
$$
\|I_{G_{r_{n}}\setminus G_{r_{n+1}}}Du
\|_{\bL_{p}(\tau)}^{p}
\leq N\big(
\|I_{G_{r_{n-1}}\setminus G_{r_{n+2}}}f^{0}
\|_{\bL_{p}(\tau)}^{p/2}
\|I_{G_{r_{n-1}}\setminus G_{r_{n+2}}}u
\|_{\bL_{p}(\tau)}^{p/2}
$$
$$
+\sum_{i=1}^{d}\|I_{G_{r_{n-1}}\setminus G_{r_{n+2}}}f^{i}
\|_{\bL_{p}(\tau)}^{p}
+\|I_{G_{r_{n-1}}\setminus G_{r_{n+2}}}g
\|_{\bL_{p}(\tau)}^{p}\big)
$$
$$
+Nr_{n+2}^{-p}\|I_{G_{r_{n-1}}\setminus G_{r_{n+2}}}u
\|_{\bL_{p}(\tau)}^{p}.
$$
  Young's inequality yields that for any constant $\chi>0$
$$
\|I_{G_{r_{n}}\setminus G_{r_{n+1}}}Du
\|_{\bL_{p}(\tau)}^{p}
\leq  Nr_{n+2}^{-p}(1+\chi) \|I_{G_{r_{n-1}}\setminus G_{r_{n+2}}}u
\|_{\bL_{p}(\tau)}^{p}
$$
$$
+N\big(r_{n+2}^{p}\chi^{-1}
\|I_{G_{r_{n-1}}\setminus G_{r_{n+2}}}f^{0}
\|_{\bL_{p}(\tau)}^{p}
 +\sum_{i=1}^{d}\|I_{G_{r_{n-1}}\setminus G_{r_{n+2}}}f^{i}
\|_{\bL_{p}(\tau)}^{p}
$$
$$
+\|I_{G_{r_{n-1}}\setminus G_{r_{n+2}}}g
\|_{\bL_{p}(\tau)}^{p}\big).
$$
We multiply both parts by $r^{\theta-d}_{n+2}$ 
and use the facts that
$r_{n-1}=2r_{n+2}$ and
on  
$G_{r_{n-1}}\setminus G_{r_{n+2}}$ the ratio $x^{1}/r_{n+2}$
satisfies
$$
1\leq x^{1}/r_{n+2}\leq 2.
$$
Then we obtain
$$
\|I_{G_{r_{n}}\setminus G_{r_{n+1}}}Du
\|_{\bL_{p,\theta}(\tau)}^{p}\leq N(1+\chi)
\|I_{G_{r_{n-1}}\setminus G_{r_{n+2}}}M^{-1}u
\|_{\bL_{p,\theta}(\tau)}^{p}
$$
$$
+  N\big(\chi^{-1}
\|I_{G_{r_{n-1}}\setminus G_{r_{n+2}}}Mf^{0}
\|_{\bL_{p,\theta}(\tau)}^{p}
 +\sum_{i=1}^{d}\|I_{G_{r_{n-1}}\setminus G_{r_{n+2}}}f^{i}
\|_{\bL_{p,\theta}(\tau)}^{p}
$$
$$
+\|I_{G_{r_{n-1}}\setminus G_{r_{n+2}}}g
\|_{\bL_{p,\theta}(\tau)}^{p}\big).
$$

Upon summing up these inequalities over $n\geq0$ we conclude
$$
\|I_{G_{r }}  D u
\|_{\bL_{p,\theta}(\tau)}^{p} 
\leq N(1+\chi)
\|I_{G_{r-1}}  M^{-1} u
\|_{\bL_{p,\theta}(\tau)}^{p}
$$
$$
+ N\big(\chi^{-1}
\|I_{G_{r_{-1} } }Mf^{0}
\|_{\bL_{p,\theta}(\tau)}^{p}
+\sum_{i=1}^{d}\|I_{G_{r_{-1} } }f^{i}
\|_{\bL_{p,\theta}(\tau)}^{p}
+ \|I_{I_{G_{r _{-1}} }}g
\|_{\bL_{p,\theta}(\tau)}^{p}\big),
$$
which after minimizing with respect
to $\chi>0$ leads to a result which is even
 somewhat sharper than \eqref{1.1.4}.
The theorem is proved.

By letting $r\to\infty$ in \eqref{1.1.4}
 we get the following.
\begin{corollary}
                                 \label{corollary 12.29.1}
If the assumptions of Theorem \ref{theorem 12.28.1}
hold with $R=\infty$, then
$$
\| D u\|_{\bL_{p,\theta}(\tau)}
\leq
N\|Mf^{0}\|_{\bL_{p,\theta}(\tau)}^{1/2}
\|M^{-1}u\|_{\bL_{p,\theta}(\tau)}^{1/2}
$$
$$
+ N\sum_{i=1}^{d}\|f^{i}\|_{\bL_{p,\theta}(\tau)}
+N\|g\|_{\bL_{p,\theta}(\tau)}
+N\|M^{-1}u\|_{\bL_{p,\theta}(\tau)},
$$
where   $N=N(d,p,\delta ,\varepsilon, \varepsilon_{1},K)$.
In particular, if $\|M^{-1}u\|_{\bL_{p,\theta}(\tau)}<\infty$,
then $u\in\frW^{1}_{p,\theta,0}(\tau)$.

\end{corollary}

Corollary \ref{corollary 12.29.1} reduces obtaining
an estimate for $\|M^{-1}u\|_{\bW^{1}_{p,\theta}(\tau)}$
to estimating $\|M^{-1}u\|_{\bL_{p,\theta}(\tau)}$.
Estimating the latter will be done by using the following
 ``energy" estimate.
Recall that
 $$
\gamma=\theta-d-p+1\quad(<0).
$$

\begin{lemma}
                                     \label{lemma 12.3.1}
Let  
$u\in\frW^{1}_{p,\theta,0}(\tau) $, 
$Mf^{0}\in \bL_{p,\theta }(\tau)$,
$ f^{i}\in \bL_{p,\theta }(\tau)$, $i=1,...,d$,
  $g=(g^{k})\in\bL_{p,\theta }(\tau)$ and assume that
\eqref{12.3.1} holds
for $t\leq\tau$ in the sense of generalized functions
on $\bR^{d}_{+}$.
Then  
$$
E\int_{0}^{\tau}
\big( \int_{\bR^{d}}\big[pM^{\gamma+1}
|u_{t}|^{p-2}u_{t}f^{0}_{t}
-p(p-1)M^{\gamma+1}|u_{t}|^{p-2}f^{i}_{t}D_{i}u_{t}
$$
$$
-p(\gamma+1)M^{\gamma }|u_{t}|^{p-2}u_{t}f^{1}_{t}
+(1/2)p(p-1)M^{\gamma+1}|u_{t}|^{p-2}|g_{t}|_{\ell_{2}}^{2}
\big]\,dx\big)\,dt 
$$
\begin{equation}
                                       \label{12.3.2}
\geq EI_{\tau<\infty}\int_{\bR^{d}}M^{\gamma+1}
|u_{\tau}| ^{p} \,dx
\end{equation}
with an equality in place of the inequality
if $\tau$ is bounded.

\end{lemma}

Proof. First of all observe that 
the concavity of the function $\log t$ implies that
$$
a_{1}^{\alpha_{1}}...a_{n}^{\alpha_{n}}\leq \alpha_{1}a_{1}
+...+\alpha_{n}a_{n}
$$
if $a_{i},\alpha_{i}\geq0$ and $\alpha_{1}+...+\alpha_{n}=1$.
It follows that for any $\kappa>0$
$$
M^{\gamma+1}
|u_{t}|^{p-1} |f^{0}_{t}|\leq \kappa M^{\theta-d}|M^{-1}u_{t}|^{p}
+NM^{\theta-d}|Mf^{0}_{t}|^{p},
$$
$$
M^{\gamma+1}|u_{t}|^{p-2}|f^{i}_{t}D_{i}u_{t}|
\leq\kappa M^{\theta-d}(|M^{-1}u_{t}|^{p}
+ |D u_{t}|^{p})+N\sum_{i=1}^{d}
 M^{\theta-d}|f^{i}_{t}|^{p},
$$
$$
M^{\gamma }|u_{t}|^{p-1}|f^{1}_{t}|
\leq \kappa M^{\theta-d}|M^{-1}u_{t}|^{p}
+NM^{\theta-d}| f^{1}_{t}|^{p},
$$
\begin{equation}
                                                \label{5.29.1}
M^{\gamma+1}|u_{t}|^{p-2}|g_{t}|_{\ell_{2}}^{2}
\leq \kappa M^{\theta-d}|M^{-1}u_{t}|^{p}
+NM^{\theta-d}|g_{t}|^{p}_{\ell_{2}} ,
\end{equation}
where the constants $N$ depend only on $\kappa$ and $p$.
The right-hand sides in these estimates are summable
over $\opar0,\tau\cbrk\times\bR^{d}$, implying that
the expectation in \eqref{12.3.2} makes perfect sense.

Next take a
nonnegative function $\phi$
of
one variable $x^{1}$ of class
$ C^{\infty}_{0}(\bR _{+})$ and notice that
$$
d(M^{(\gamma+1)/p}u_{t}\phi)=
 \big(M^{(\gamma+1)/p}f^{0}_{t}\phi-M^{(\gamma+1)/p}
f^{1}_{t}\phi'
$$
$$
-(\gamma+1)p^{-1}
M^{(\gamma+1)/p-1}\phi f^{1}_{t}
+D_{i}(
M^{(\gamma+1)/p}f^{i}_{t}\phi)\big)\,dt
+M^{(\gamma+1)/p}g^{k}_{t}\phi\,dw^{k}_{t}.
$$
This equation holds in $\bR^{d}$ rather than only in $\bR^{d}_{+}$.
Hence, by Corollary 2.2 of~\cite{Ito} 
$$
E\int_{0}^{\tau}\big(\int_{\bR^{d}}M^{\gamma+1}
\big\{p|u_{t}|^{p-2}u_{t}\phi^{p-1}[\phi f^{0}_{t}-
pf^{1}_{t}\phi'-(\gamma+1) 
 M^{-1}\phi f^{1}_{t}]
$$
$$
-p(p-1) |u_{t}|^{p-2}\phi^{p-1}f^{i}_{t}
 \phi  D_{i}u_{t} 
+(1/2)p(p-1) \phi^{p}|u_{t}|^{p-2}|g_{t}|_{\ell_{2}}
^{2}\big\}\,dx\big)\,dt
$$
\begin{equation}
                                         \label{6.11.1}
\geq EI_{\tau<\infty}\int_{\bR^{d}}M^{\gamma+1}
|u_{\tau}| ^{p}\phi^{p}\,dx
\end{equation}
with an equality in place of the inequality
if $\tau$ is bounded.

By recalling what was said in the beginning of the proof
and having in mind the dominated convergence
theorem and Fatou's lemma
we easily see that, to prove inequality \eqref{12.3.2}, now
 it suffices to find a sequence of
$\phi_{n}\in C^{\infty}_{0}(\bR _{+})$ such that
$0\leq\phi_{n}\leq 1$, $\phi_{n}\to1$, and
$$
E\int_{0}^{\tau} \int_{\bR^{d}}M^{\gamma+1}
 |u_{t}|^{p-1} | f^{1}_{t}\phi_{n}'|\,dxdt\to0.
$$
Furthermore, since estimates \eqref{5.29.1} imply that
$$
E\int_{0}^{\tau} \int_{\bR^{d}}M^{\gamma }
 |u_{t}|^{p-1} | f^{1}_{t} |\,dxdt<\infty,
$$
the dominated convergence theorem shows that
it suffices to find a sequence of
$\phi_{n}\in C^{\infty}_{0}(\bR _{+})$ such that
$0\leq\phi_{n}\leq 1$, $\phi_{n}\to1$,
$M\phi_{n}'$ are uniformly bounded, and
$M\phi_{n}'\to0$ in $\bR_{+}$.

To construct such a sequence, take some nonnegative
$\eta,\zeta\in C^{\infty}_{0}
(\bR)$ such that $\eta=0$ near the origin, $\eta(x)=1$ for
$x\geq1$,  $\zeta=1$ near the origin, and $\eta,\zeta\leq1$.
 Then define
$\phi_{n}(x)=\eta(nx)\zeta( x/n)$.
The reader will easily check that the required 
properties are satisfied.

To prove that \eqref{12.3.2} holds with the equality
sign if $\tau$ is bounded, we write
\eqref{6.11.1} with the equality sign and
pass to the limit by the dominated convergence theorem
knowing already that the right-hand side of 
\eqref{12.3.2} is finite.
The lemma is proved.

\begin{corollary}
                                    \label{corollary 5.29.1}
Let  Assumptions \ref{assumption 3.6.1}
and \ref{assumption 3.6.2}  be satisfied. Let
$u\in\frW^{1}_{p,\theta,0}(\tau) $, 
$Mf^{0}\in \bL_{p,\theta }(\tau)$,
$ f^{i}\in \bL_{p,\theta }(\tau)$, $i=1,...,d$,
  $g=(g^{k})\in\bL_{p,\theta }(\tau)$ and assume that
$u$ satisfies \eqref{11.13.1} for $t\leq\tau$. Then
for any constant $\chi>0$ there exist   constants
$N^{*}=N^{*}(d,p,\delta)$ and
$N=N(\chi,d,p,\delta  )$ such that
$$
 p(p-1)E\int_{0}^{\tau}
  \int_{\bR^{d}}M^{\gamma+1}|u_{t}|^{p-2}(a^{ij}_{t}
-\alpha^{ij}_{t})(D_{i}u_{t})D_{j}u_{t}
 \,dxdt
$$
$$
+p(\gamma+1)E\int_{0}^{\tau}
  \int_{\bR^{d}}M^{\gamma }|u_{t}|^{p-2}u_{t}a^{i1}_{t}D_{i}u_{t}
 \,dxdt
$$
$$
\leq N(\|Mf^{0}\|_{\bL_{p,\theta}(\tau)}^{p}
+\sum_{i=1}^{d}\| f^{i}\|_{\bL_{p,\theta}(\tau)}^{p}
+(1+K^{p})\|g\|_{\bL_{p,\theta}(\tau)}^{p})
$$
\begin{equation}
                                                  \label{5.30.4}
+ [ N^{*}K(1+K)+\chi] I,
\end{equation}
where
$$
I=E\int_{0}^{\tau}
  \int_{\bR^{d}}(M^{\gamma-1}|u_{t}|^{p}+M^{\gamma+1}
|u_{t}|^{p-2}|Du_{t}|^{2})\,dxdt\leq N^{*}\|M^{-1}u\|^{p}
_{\bW^{1}_{p,\theta}(\tau)}.
$$
\end{corollary}

To derive this result observe that by Lemma \ref{lemma 12.3.1}
$$
E\int_{0}^{\tau}
\big( \int_{\bR^{d}}\big[pM^{\gamma+1}
|u_{t}|^{p-2}u_{t}(b^{i}_{t}D_{i}u_{t}+c_{t}u_{t}+f^{0}_{t})
$$
$$
-p(p-1)M^{\gamma+1}|u_{t}|^{p-2}(
a^{ij}_{t}D_{i}u_{t}+a^{j}_{t}u_{t}+f^{j}_{t})D_{j}u_{t}
$$
$$
-p(\gamma+1)M^{\gamma }|u_{t}|^{p-2}u_{t}
(a^{i1}_{t}D_{i}u_{t}+a^{1}_{t}u_{t}+f^{1}_{t})
$$
$$
+(1/2)p(p-1)M^{\gamma+1}|u_{t}|^{p-2}|
\sigma^{i\cdot}_{t}D_{i}u_{t}+\nu_{t}u_{t}+g_{t}|_{\ell_{2}}^{2}
\big]\,dx\big)\,dt\geq0,
$$
which after rearranging the terms becomes
$$
 p(p-1)E\int_{0}^{\tau}
  \int_{\bR^{d}}M^{\gamma+1}|u_{t}|^{p-2}(a^{ij}_{t}
-\alpha^{ij}_{t})(D_{i}u_{t})D_{j}u_{t}
 \,dxdt
$$
$$
+p(\gamma+1)E\int_{0}^{\tau}
  \int_{\bR^{d}}M^{\gamma }|u_{t}|^{p-2}u_{t}a^{i1}_{t}D_{i}u_{t}
 \,dxdt
$$
$$
\leq E\int_{0}^{\tau}
  \int_{\bR^{d}}\big[M^{\gamma }|u_{t}|^{p-2}u_{t}
A^{i}_{t}D_{i}u_{t} +M^{\gamma-1}|u_{t}|^{p} B_{t}
\big]\,dxdt
$$
$$
+E\int_{0}^{\tau}
  \int_{\bR^{d}}\big[M^{\gamma }|u_{t}|^{p-2}u_{t} F_{t}
+M^{\gamma+1}|u_{t}|^{p-2}G_{t}+
M^{\gamma+1}|u_{t}|^{p-2}H^{i}_{t}D_{i}u_{t}\big]\,dxdt,
$$
where
$$
A^{i}_{t}=pMb^{i}_{t}-p(p-1)Ma^{i}_{t}+p(p-1)(\sigma^{i\cdot}
_{t},M\nu_{t})_{\ell_{2}},
$$
 $$
B_{t}=
 pM^{2}c_{t}-p(\gamma+1)Ma^{1}_{t}
+(1/2)p(p-1)M^{2}|\nu_{t}|_{\ell_{2}}^{2}  ,
$$
$$
F_{t}=pM f^{0}_{t} -p(\gamma+1)f^{1}_{t}
+p(p-1)(M\nu_{t},g_{t})_{\ell_{2}} ,
$$
$$
G_{t}=(1/2)p(p-1) | g_{t}|^{2}
_{\ell_{2}} ,
$$
$$
H^{i}_{t}=  p(p-1)(\sigma^{i\cdot}_{t},g_{t})
_{\ell_{2}}-p(p-1)f^{i}_{t}.
$$
To estimate the first expectation on the right,
one uses the following simple estimates
$$
|A^{i}_{t}|\leq N^{*}K,\quad|B_{t}|\leq N^{*}K(1+K),
$$
$$
M^{\gamma}|u_{t}|^{p-1}|Du_{t}|
=(M^{(\gamma-1)/2}|u_{t}|^{p/2})(|Du_{t}|M^{(\gamma+1)/2}
|u_{t}|^{(p-2)/2})
$$
$$
\leq M^{\gamma-1}|u_{t}|^{p}+M^{\gamma+1}
|u_{t}|^{p-2}|Du_{t}|^{2}=M^{\theta-d}|M^{-1}u_{t}|^{p}
$$
$$
+M^{\theta-d}|M^{-1}u_{t}|^{p-2}|Du_{t}|^{2}
\leq 2M^{\theta-d}|M^{-1}u_{t}|^{p} 
+M^{\theta-d}|Du_{t}|^{p}.
$$
The second expectation is estimated by using inequalities
like \eqref{5.29.1}. For instance,
$$
M^{\gamma+1}|u_{t}|^{p-2}|D u_{t}|\,|H _{t}|
$$
$$
=(M^{(\gamma-1)(p-2)/(2p)}|u_{t}|^{(p-2)/2})
(M^{(\gamma+1)/2}|u_{t}|^{(p-2)/2}|Du_{t}|)
(M^{(\theta-d)/p}|H_{t}|
$$
$$
\leq\chi(M^{\gamma-1}|u_{t}|^{p}+M^{\gamma+1}
|u_{t}|^{p-2}|Du_{t}|^{2})+NM^{\theta-d}|H_{t}|^{p}.
$$

Now we prepare to estimate from below
the left-hand side of \eqref{5.30.4} in terms
of a quantity equivalent to
 $\|M^{-1}u\|_{\bL _{p,\theta}(\tau)}$.
The following two results will not be used in the
  proof of Theorem \ref{theorem 5.30.1}.

\begin{lemma}
                                          \label{lemma 6.1.2}
Let   $\beta,\varepsilon\in(0,\infty)$
  be some constants
and let $\bar{a}$ be a measurable bounded $\bR^{d}$-valued
function on $\bR^{d}_{+}$ such that  
\begin{equation}
                                                \label{6.1.4}
|\bar{a}(x)-\bar{a}(y)|\leq\beta
\end{equation}
whenever $x,y\in\bR^{d}_{+}$ and $|x-y|\leq 
 \varepsilon(x^{1}\wedge y^{1})$.
Then for any  $u\in M W^{1}_{p,\theta}$
  we have
\begin{equation}
                                                \label{6.2.3}
\big|I+p^{-1}\gamma 
\int_{\bR^{d}_{+}} 
\bar{a}^{1} M^{\gamma-1}
|u|^{p}\,dx \big|
\leq N\beta\|M^{-1}u\|_{W^{1}_{p,\theta}}^{p},
\end{equation}
where $N=N(d,p,\theta,\varepsilon )$  and
$$
I:=  \int_{\bR^{d}}M^{\gamma }|u |^{p-2}u \bar{a}^{i}D_{i} u
 \,dx. 
$$

\end{lemma}

Proof. Since $C^{\infty}_{0}(\bR^{d}_{+})$ is dense
in $MW^{1}_{p,\theta}$ we may assume that $u\in
C^{\infty}_{0}(\bR^{d}_{+})$.
Take a nonnegative
$\zeta\in C^{\infty}_{0}(\bR^{d}_{+})$ with unit integral
and such that $\zeta(x)=0$ if $x^{1}
\not\in(1,1+\varepsilon/2)$
or $|x'|\geq\varepsilon/2 $.
For   $y \in\bR^{d }_{+}$
define
$$
\zeta^{y}(x)=(x^{1})^{\gamma+1}
\zeta(y^{1}x^{1},y^{1}(y'-x'))(y^{1})^{d-1}.
$$
Observe that for $x\in\bR^{d}_{+}$
\begin{equation}
                                            \label{6.1.5} 
\int_{\bR^{d}_{+}}  \zeta^{y}(x) \,dy=
(x^{1})^{\gamma }\int_{\bR^{d}_{+}}  \zeta (y) \,dy=
(x^{1})^{\gamma }.
\end{equation}

It follows that
$$
I
=\int_{\bR^{d}_{+}} I(y)\,dy,
$$
where $I(y)=I_{1}(y)+I_{2}(y)$,
$$
I_{1}(y)=\int_{\bR^{d}_{+}} \zeta^{y} 
|u |^{p-2}u \bar{a}^{i }(\bar{y}) D_{i }u  \,dx,
\quad\bar{y}=((y^{1})^{-1},y'),
$$
$$
I_{2}(y)=\int_{\bR^{d}_{+}} \zeta^{y} 
|u |^{p-2}u [\bar{a}^{i }(x)-
\bar{a}^{i }(\bar{y})] D_{i }u \,dx.
$$

By the choice of $\zeta$ we have that if $\zeta^{y}(x)\ne0$,
then $1<y^{1}x^{1}<1+\varepsilon/2$ and
$y^{1}|y'-x'|<\varepsilon/2 $ implying that
$$
 \bar{y}^{1}< x^{1}<(1+\varepsilon/2)\bar{y}^{1},\quad
|\bar{y}'-x'|<\bar{y}^{1}\varepsilon/2 
=(\varepsilon/2)(x^{1}\wedge \bar{y}^{1}),
$$
$$
  0<x^{1}-\bar{y}^{1}<\bar{y}^{1}\varepsilon/2,\quad 
|x^{1}-\bar{y}^{1}|<(\varepsilon/2)(x^{1}\wedge \bar{y}^{1}),
$$
$$
|x-\bar{y}|< \varepsilon (x^{1}\wedge \bar{y}^{1}),\quad
| \bar{a} (x)-
\bar{a} (\bar{y})|\leq \beta.
$$

Hence,
$$
|I_{2}(y)|\leq  \beta
\int_{\bR^{d}_{+}} \zeta^{y} 
|u |^{p-1}|D u | \,dx,
$$
$$
\int_{\bR^{d}_{+}}|I_{2}(y)|\,dy
\leq  \beta\int_{\bR^{d}_{+}}M^{\gamma}|u |^{p-1}|D u | \,dx
\leq N\beta\|M^{-1}u\|_{W^{1}_{p,\theta}}^{p},
$$
\begin{equation}
                                                \label{6.2.1}
\big|I-\int_{\bR^{d}_{+}} I_{1}(y) \,dy\big|
\leq N\beta\|M^{-1}u\|_{W^{1}_{p,\theta}}^{p}.
\end{equation}

To deal with $I_{1}(y)$ we integrate by parts observing that
$$
|u|^{p-2}uD_{i}u=p^{-1}D_{i}(|u|^{p}).
$$
Then we find
$$
I_{1}(y)=-p^{-1}\int_{\bR^{d}_{+}}(D_{i}\zeta^{y})
\bar{a}^{i}(\bar{y})
|u|^{p}\,dx=-p^{-1}J_{1}(y)-p^{-1}J_{2}(y),
$$
where
$$
J_{1}(y)=\int_{\bR^{d}_{+}}(D_{i}\zeta^{y})[\bar{a}^{i}(\bar{y})-
\bar{a}^{i}(x)]
|u|^{p}\,dx,
$$
$$
J_{2}(y)=\int_{\bR^{d}_{+}}(D_{i}\zeta^{y}) 
\bar{a}^{i} 
|u|^{p}\,dx.
$$

As is easy to see
$$
\int_{\bR^{d}_{+}}D_{i}\zeta^{y}\,dy=D_{i}((x^{1})^{\gamma})
=\delta^{i1}\gamma (x^{1})^{\gamma-1},
$$
$$
\int_{\bR^{d}_{+}}J_{2}(y)\,dy=
\gamma \int_{\bR^{d}_{+}} 
\bar{a}^{1} M^{\gamma-1}
|u|^{p}\,dx 
$$
and by \eqref{6.2.1}
\begin{equation}
                                                \label{6.2.2}
\big|I+p^{-1}\gamma 
\int_{\bR^{d}_{+}} 
\bar{a}^{1} M^{\gamma-1}
|u|^{p}\,dx \big|
\leq N\beta\|M^{-1}u\|_{W^{1}_{p,\theta}}^{p}
+p^{-1}\int_{\bR^{d}_{+}} |J_{1}(y)| \,dy.
\end{equation}

Furthermore,
$$
|J_{1}(y)|\leq\beta
\int_{\bR^{d}_{+}}|D \zeta^{y}|\,
|u|^{p}\,dx.
$$
Here
$$
|D\zeta^{y}(x)|\leq|\gamma+1|(x^{1})^{-1}\zeta^{y}(x)
+ \hat{\zeta}^{y}(x)y^{1} ,
$$
$$
\hat{\zeta}^{y}(x):=(x^{1})^{\gamma+1}|D\zeta|
(y^{1}x^{1},y^{1}(y'-x'))(y^{1})^{d-1},
$$
$$
\int_{\bR^{d}_{+}}|D\zeta^{y}(x)|\,dy
\leq|\gamma+1|(x^{1})^{\gamma-1}
+(x^{1})^{\gamma-1}\int_{\bR^{d}_{+}}|D\zeta (y)|y^{1}\,dy
=N(x^{1})^{\gamma-1}.
$$
 It follows that
$$
\int_{\bR^{d}_{+}}|J_{1}(y)|\,dy\leq N\beta
\int_{\bR^{d}_{+}}M^{\gamma-1}\,
|u|^{p}\,dx=N\beta\|M^{-1}u\|_{L_{p,\theta}},
$$
which after being combined with \eqref{6.2.2}
leads to \eqref{6.2.3} and proves the lemma.

The following lemma  is a simple
consequence of Lemma 6.6 of \cite{Kr08},
where the estimate is stronger. The proof
of Lemma 6.6 of \cite{Kr08} follows the same lines as that
of Lemma \ref{lemma 6.1.2}. Lemma \ref{lemma 1.4.02}
will be used  for $\bar{a}^{ij}=(a^{11}_{t})^{-1}
\hat{a}^{i1}_{t}\hat{a}^{j1}_{t}$.

\begin{lemma}
                                         \label{lemma 1.4.02}
Let     $\beta,\varepsilon\in(0,\infty)$
  be some constants and let
  $\bar{a}(x)$ be a measurable 
function given on $\bR^{d}_{+}$ with values in the set
of symmetric nonnegative matrices and such that
$|\bar{a}^{ij}|\leq \delta^{-1}$
and 
\begin{equation}
                                             \label{2.23.1}
|\bar{a}^{ij}(x)-\bar{a}^{ij}(y)|\leq\beta
\end{equation}
whenever $x,y\in\bR^{d}_{+}$ and $|x-y|\leq 
 \varepsilon(x^{1}\wedge y^{1})$.
Then for any  $u\in M W^{1}_{p,\theta}$ and $
\chi>0$ and $\kappa\in(0,1]$
  we have
$$
 \int_{\bR^{d}_{+}}M^{\gamma+1}
|u |^{p-2}\bar{a}^{ij}  (D_{i}u )D_{j}u  \,dx
$$
$$
\geq (1-\kappa)\gamma^{2}p^{-2}
\int_{\bR^{d}_{+}}M^{\gamma-1}\bar{a}^{11} 
|u|^{p}\,dx
$$
\begin{equation}
                                                   \label{1.11.2}
-N \big((\varepsilon^{-1}R+1)
\beta  +\kappa^{-1}\chi\big)
\|M^{-1}u\|_{W^{1}_{p,\theta}}^{p},
\end{equation}
where  $N=N(d,p,\delta,\theta)$ and
 $\ln R=N(d,p)\chi^{-1/2}$.

\end{lemma}

\mysection{Proof of Theorems \protect\ref{theorem 5.30.1}
and \protect\ref{theorem 6.2.1}}

                                             \label{section 6.4.2}

With start with a theorem that says that to prove the solvability of
\eqref{11.13.1} we only need to have an a priori
estimate of the lowest norm of $u$.

\begin{theorem}
                                    \label{theorem 1.30.1}
Let   Assumptions \ref{assumption 3.6.1}   
  through \ref{assumption 8.10.1} 
be satisfied.
 Assume that there
is a constant $N_{0}<\infty$ such that for any 
$\lambda\in[0,1]$, $u\in\frW^{1}_{p,\theta,0}(\tau)$,
and $f^{0}, ..., f^{d}$  and
$g=(g^{k})$, satisfying
\begin{equation}
                                               \label{1.30.6} 
  Mf^{0}, f^{i},
g=(g^{k})\in\bL_{p,\theta}(\tau),\quad i=1,...,d, 
\end{equation}
we have the a priori estimate
\begin{equation}
                                            \label{1.30.3}
\| M^{-1}u\|_{\bL_{p,\theta}(\tau)}\leq
N_{0}\big(\|Mf^{0}\|_{\bL_{p,\theta}(\tau)}
+\sum_{i=1}^{d}\|f^{i}\|_{\bL_{p,\theta}(\tau)}
+\|g\|_{\bL_{p,\theta}(\tau)}\big)
\end{equation}
provided that
$$
du_{t}=(\lambda\Lambda^{k}_{t}
u_{t}+g^{k}_{t})\,dw^{k}_{t}
$$
\begin{equation}
                                            \label{1.30.4}
+ 
 [(\lambda L_{t}+(1-\lambda)\Delta)u_{t}
+f^{0}_{t}+D_{i}f^{i}_{t} ]\,dt ,\quad t\leq\tau ,
\end{equation}
in $\bR^{d}_{+}$ (estimate \eqref{1.30.3} is {\em not\/}
supposed to hold if there is no solution $u
\in\frW^{1}_{p,\theta,0}(\tau)$ of \eqref{1.30.4}).

Then  for any  $f^{0}, ..., f^{d}$, and
$g=(g^{k})$ satisfying \eqref{1.30.6}
there exists a unique
$u
\in\frW^{1}_{p,\theta,0}(\tau)$ satisfying \eqref{11.13.1}
in $\bR^{d}_{+}$ for $t\leq\tau$.  Furthermore, for this solution
\begin{equation}
                                               \label{1.30.5}
\| Du\|_{\bL_{p,\theta}(\tau)}  
\leq N\big(\| Mf^{0}\|_{\bL_{p,\theta}(\tau)}
+\sum_{i=1}^{d}\|  f^{i}\|_{\bL_{p,\theta}(\tau)}
+ \| g\|_{\bL_{p,\theta}(\tau)}\big),
\end{equation}
where $N$ depends only on
$d,p,\delta , K,\varepsilon,\varepsilon_{1} $, and $N_{0}$.
\end{theorem}

Proof. We call a $\lambda\in[0,1]$ ``good"
if  for any  
for any  $f^{0}, ..., f^{d}$, and
$g=(g^{k})$ satisfying
\eqref{1.30.6}
there exists a unique
$u
\in\frW^{1}_{p,\theta,0}(\tau)$ satisfying \eqref{1.30.4}
in $\bR^{d}_{+}$. By Corollary \ref{corollary 12.29.1}
and   assumption \eqref{1.30.3}   
estimate \eqref{1.30.5} holds for solutions of 
\eqref{1.30.4} if $\lambda$ is a ``good" point.
It follows that to prove the theorem it suffices
to prove that all points of $[0,1]$ are ``good".

We are going to use the method of continuity
observing that the fact that
the point $0$ is ``good" is known
from \cite{KL}
(or is easily obtained as suggested after
 \eqref{6.4.1}). We will achieve our goal if we show that
there exists a constant $\mu>0$ such that
if $\lambda_{0}$ is a ``good" point, then
all points in the interval $[\lambda_{0}-\mu, 
\lambda_{0}+\mu]\cap[0,1]$ are ``good". So fix a ``good"
point $\lambda_{0}$ and fix some
$f^{0}, ..., f^{d}$, and
$g=(g^{k})$ satisfying
\eqref{1.30.6}.

For any $v\in M\bW^{1}_{p,\theta}(\tau)$
consider the equation
$$
du_{t}=[(\lambda_{0} L_{t}+(1-\lambda_{0})\Delta)u_{t}
+(\lambda-\lambda_{0})(L_{t}-\Delta)v_{t}
+D_{i}f^{i}_{t}+f^{0}_{t})\,dt
$$
\begin{equation}
                                          \label{12.7.4}
+(\lambda_{0} \Lambda^{k}_{t}u_{t}
+(\lambda-\lambda_{0})\Lambda^{k}v_{t}+g^{k}_{t})\,dw^{k}_{t}.
\end{equation}
Observe that
$$
(L_{t}-\Delta)v_{t}=D_{j}\big((a^{ij}-\delta^{ij})
D_{i}v_{t}+a^{j}_{t}v_{t}\big)+b^{i}_{t}D_{i}v_{t}
+cv_{t},
$$
where by assumption
$$
|(a^{ij}-\delta^{ij})
D_{i}v_{t}|\leq N|Dv_{t}|,\quad |a^{j}_{t}v_{t}|
\leq NM^{-1}|v_{t}|, \quad M|b^{i}_{t}D_{i}v_{t}|
\leq N|Dv_{t}|,
$$
$$
  M|cv_{t}|\leq NM^{-1}|v_{t}|,\quad
|\Lambda^{\cdot}v_{t}|_{\ell_{2}}\leq N(|Dv_{t}|+M^{-1}|v_{t}|)
$$
and the right-hand sides in these estimates are 
in $\bL_{p,\theta}(\tau)$.
Hence by the assumption that $\lambda_{0}$
is ``good", equation \eqref{12.7.4} has a unique solution
$u\in\frW^{1}_{p,\theta,0}(\tau)$ 
($\subset M\bW^{1}_{p,\theta}(\tau)$).

In this way, for $f^{j}$ and $g$ being fixed,
 we define a mapping
$v\to u$ in the space $M\bW^{1}_{p,\theta}(\tau)$. 
It is important to keep in mind
 that the image $u$ of 
$v\in M\bW^{1}_{p,\theta}(\tau)$ is always in 
$\frW^{1}_{p,\theta,0}(\tau)$.
Take
$v',v''\in M\bW^{1}_{p,\theta}(\tau)$ and let $u',u''$ be their 
corresponding images. Then $u:=u'-u''$ satisfies
$$
du_{t}=[(\lambda_{0} L_{t}+(1-\lambda_{0})\Delta)u_{t}
+(\lambda-\lambda_{0})(L_{t}-\Delta)v_{t})\,dt
$$
$$
+(\lambda\Lambda^{k}_{t}u_{t}
+(\lambda-\lambda_{0})\Lambda^{k}v_{t} )\,dw^{k}_{t},
$$
where $v=v'-v''$. It follows by \eqref{1.30.3} and \eqref{1.30.5}
that
$$
\|M^{-1}u\|_{\bW^{1}_{p,\theta}(\tau)}\leq N|\lambda-\lambda_{0}|
\,\|M^{-1}v\|_{\bW^{1}_{p,\theta}(\tau)}
$$
with a constant $N$ independent of $f$, $g$, $v'$, $v''$, $\lambda_{0}$,
and $\lambda$. 
For $\lambda$ sufficiently close
to $\lambda_{0}$, our mapping is a contraction
and, since $M\bW^{1}_{p}(\tau)$ is a Banach
space, 
the mapping
  has a fixed point. This fixed point is in $\frW^{1}_{p,\theta,0}
(\tau)$ and,
obviously, satisfies \eqref{1.30.4}.
 As is explained above, this proves the theorem.

{\bf Proof of Theorem \ref{theorem 5.30.1}}.
 According to Theorem \ref{theorem 1.30.1}  it suffices to find
$K=K(d,p,\delta,
\bar{\delta},\theta,\varepsilon,\varepsilon_{1})>0$ such that
Assumptions \ref{assumption 3.6.1} through 
 \ref{assumption 8.10.1}
would imply that  
 \eqref{1.30.3} holds 
for any solution $u\in\frW^{1}_{p,\theta,0}(\tau)$
of \eqref{11.13.1} for $t\leq\tau$ and
  $N_{0}$
depends only on $d,p,\delta,\theta,\bar{\delta}$, $\varepsilon $, and 
$\varepsilon_{1}$.  From the start we will only consider
$K\leq1$. This assumption allows us to eliminate $K$
from the lists of what $N$'s depend on in Theorem \ref{theorem 1.30.1}
and Corollary \ref{corollary 12.29.1}.

 By H\"older's inequality
$$
I:=\big|E\int_{0}^{\tau}
  \int_{\bR^{d}}M^{\gamma }|u_{t}|^{p-2}u_{t}a^{i1}_{t}D_{i}u_{t}
 \,dxdt\big|\leq I_{1}^{1/2}I_{2}^{1/2},
$$
where
$$
I_{1}=E\int_{0}^{\tau}
  \int_{\bR^{d}}M^{\gamma+1}|u_{t}|^{p-2}\big(\sum_{i}
a^{i1}_{t}D_{i}u_{t}\big)^{2}
 \,dxdt,
$$
and $I_{2}=\|M^{-1}u\|^{p}_{\bL_{p,\theta}(\tau)}$. By assumption
\eqref{6.1.1}
$$
I_{1}\leq\bar{\delta}^{-1}
E\int_{0}^{\tau}
  \int_{\bR^{d}}M^{\gamma+1}|u_{t}|^{p-2} 
(a^{ij}_{t}-\alpha^{ij}_{t})(D_{j}u_{t})D_{i}u_{t} 
 \,dxdt=:\bar{\delta}^{-1}I_{3}.
$$
By Lemma 6.1 of \cite{Kr8} (Hardy's inequality)
and
Assumption \ref{assumption 3.6.2} we have
\begin{equation}
                                             \label{6.1.2}
\gamma^{2}I_{2}\leq p^{2}
E\int_{0}^{\tau}
  \int_{\bR^{d}}M^{\gamma+1}|u_{t}|^{p-2}|
 D u_{t}|^{2}
 \,dxdt\leq p^{2}\delta ^{-1}I_{3}.
\end{equation}

It follows that 
$$
I\leq \bar{\delta}^{-1/2}\delta^{-1/2}p|\gamma|^{-1}I_{3},
$$
so that
the left hand side of \eqref{5.30.4}
dominates
$$
p(p-1)I_{3}-p|\gamma+1|\bar{\delta}^{-1/2}|\gamma|
^{-1}p\delta^{-1/2}I_{3}.
$$
By assumption the sum of the coefficients of $I_{3}$
is strictly positive. Therefore, a strictly positive factor
of $I_{3}$ admits an estimate in terms of the right-hand side 
of \eqref{5.30.4}. Estimate \eqref{6.1.2} shows that the same is
true for $I_{2}$. 
In other words,
$$
\|M^{-1}u\|^{p}_{\bL_{p,\theta}(\tau)}
\leq NJ
+ [ N^{*}K(1+K)+\chi] 
\|M^{-1}u\|^{p}_{\bW^{1}_{p,\theta}(\tau)},
$$
where $\chi>0$ is arbitrary, $N=N(\chi,d,p,\delta,\bar{\delta},
\theta)$, $N^{*}=N^{*}(d,p,\delta,\bar{\delta},
\theta)$ and
$$
J=\|Mf^{0}\|_{\bL_{p,\theta}(\tau)}^{p}
+\sum_{i=1}^{d}\| f^{i}\|_{\bL_{p,\theta}(\tau)}^{p}
+ \|g\|_{\bL_{p,\theta}(\tau)}^{p}.
$$
 Upon combining this with Corollary \ref{corollary 12.29.1}
we find
$$
\|M^{-1}u\|^{p}_{\bW^{1}_{p,\theta}(\tau)}
\leq NJ
+ [ N^{*}K(1+K)+\chi] 
\|M^{-1}u\|^{p}_{\bW^{1}_{p,\theta}(\tau)},
$$
 where
$N=N(\chi,d,p,\delta,\bar{\delta},
\theta,\varepsilon,\varepsilon_{1})$
 and $N^{*}=N^{*}(d,p,\delta,\bar{\delta},
\theta,\varepsilon,\varepsilon_{1})$.
 
Now it is clear how to find   $
\chi >0$ and $K>0$, depending only on
$d,p,\delta,\bar{\delta},
\theta,\varepsilon$, and $\varepsilon_{1}$, so that 
the last estimate 
 would imply that the estimate
$$
\|M^{-1}u\|^{p}_{\bL_{p,\theta}(\tau)}\leq
\|M^{-1}u\|^{p}_{\bW^{1}_{p,\theta}(\tau)}
\leq N_{0}J,
$$
  implying \eqref{1.30.3}, holds
for any solution $u\in\frW^{1}_{p,\theta,0}(\tau)$
of \eqref{11.13.1} with
  $N_{0}$
depending only on
$d,p,\delta,\bar{\delta},
\theta,\varepsilon$, and $\varepsilon_{1}$.
The theorem is proved.

{\bf Proof of Theorem \ref{theorem 6.2.1}}.
 As in the above proof, given
$d,p,\delta,\tilde{\delta},
\theta,\varepsilon$, and $\varepsilon_{1}$, it suffices
to show how to find $K>0$ and $\beta_{2}>0$ in such a way that
Assumptions \ref{assumption 3.6.1}  
through \ref{assumption 6.2.1} would allow us to derive
\eqref{1.30.3}. Again   without loss of generality we assume
that $K,\beta_{2}\leq1$.  By Lemma \ref{lemma 6.1.2}
$$
p(\gamma+1)E\int_{0}^{\tau}
  \int_{\bR^{d}}M^{\gamma }|u_{t}|^{p-2}u_{t}a^{i1}_{t}D_{i}u_{t}
 \,dxdt
$$
$$
\geq -\gamma(\gamma+1)
E\int_{0}^{\tau}
  \int_{\bR^{d}}a^{11}_{t}
M^{\gamma-1}|u_{t}|^{p } \,dxdt-N\beta_{2}
\|M^{-1}u\|_{W^{1}_{p,\theta}(\tau)}^{p},
$$
where $N=N(d,p,\theta )$.
By Assumption \ref{assumption 6.2.1} and
Lemma \ref{lemma 1.4.02} for $\bar{a}^{ij}=(a^{11}_{t})^{-1}
\hat{a}^{i1}_{t}\hat{a}^{j1}_{t}$ we have
$$
 p(p-1)E\int_{0}^{\tau}
  \int_{\bR^{d}}M^{\gamma+1}|u_{t}|^{p-2}(a^{ij}_{t}
-\alpha^{ij}_{t})(D_{i}u_{t})D_{j}u_{t}
 \,dxdt
$$
$$
\geq  p(p-1)\tilde{\delta}E\int_{0}^{\tau}
  \int_{\bR^{d}}M^{\gamma+1}|u_{t}|^{p-2}\bar{a}^{ij}_{t}
 (D_{i}u_{t})D_{j}u_{t}
 \,dxdt
$$
$$
\geq p^{-1}(p-1)\tilde{\delta}(1-\kappa)\gamma^{2}
E\int_{0}^{\tau}
  \int_{\bR^{d}}a^{11}_{t}
M^{\gamma-1}|u_{t}|^{p } \,dxdt 
$$
$$
-N \big(( R+1)
\beta_{2}  +\kappa^{-1}\chi\big)
\|M^{-1}u\|_{W^{1}_{p,\theta}(\tau)}^{p},
$$
where $N=N(d,p,\delta,\tilde{\delta},\theta)$, 
$\ln R=N(d,p)\chi^{-1/2}$,
and $\kappa\in(0,1]$ and $\chi>0$ are arbitrary.

 Observe that, as $\kappa\downarrow0$,
$$
-\gamma(\gamma+1)+p^{-1}(p-1)\tilde{\delta}(1-\kappa)\gamma^{2}
\to
-\gamma[\gamma+1 -p^{-1}(p-1)\tilde{\delta} \gamma].
$$
The latter is a strictly positive constant since $\gamma<0$
and
$$
\gamma+1 +p^{-1}(p-1)\tilde{\delta} \gamma=
\frac{p-\tilde{\delta}p+\tilde{\delta}}{p}
\big[\theta+\frac{p}{p-\tilde{\delta}p+\tilde{\delta}}-d-p+1\big]>0
$$
by Assumption \ref{assumption 6.2.1}. It follows
by \eqref{5.30.4}
that after fixing $\kappa=\kappa(d,p,\theta,\tilde{\delta})\in(0,1]$
appropriately we can find an
  $N=N(d,p,\theta,\tilde{\delta},\delta )$ such that
for any $\chi>0$
$$
\|M^{-1}u\|^{p}_{\bL_{p,\theta}(\tau)}
\leq
N \big(( R+1)
\beta_{2}  +K+\chi\big)
\|M^{-1}u\|_{W^{1}_{p,\theta}(\tau)}^{p}
$$
\begin{equation}
                                        \label{6.2.8}
+N^{*}(\|Mf^{0}\|_{\bL_{p,\theta}(\tau)}^{p}
+\sum_{i=1}^{d}\| f^{i}\|_{\bL_{p,\theta}(\tau)}^{p}
+(1+\beta^{p})\|g\|_{\bL_{p,\theta}(\tau)}^{p}),
\end{equation}
where $N^{*}=N^{*}
 (d,p,\theta,\tilde{\delta},\delta,\chi)$.

By \eqref{6.2.8} and
 Corollary \ref{corollary 12.29.1}, for any $\chi>0$,
$$
\|M^{-1}u\|_{\bW^{1}_{p,\theta}(\tau)}^{p}
\leq N_{1}  ( (R+1) 
\beta_{2}+K  +\chi\big)
\|M^{-1}u\|_{W^{1}_{p,\theta}(\tau)}^{p}
$$
\begin{equation}
                                        \label{6.2.7}
+N_{2}(\|Mf^{0}\|_{\bL_{p,\theta}(\tau)}^{p}
+\sum_{i=1}^{d}\| f^{i}\|_{\bL_{p,\theta}(\tau)}^{p}
+ \|g\|_{\bL_{p,\theta}(\tau)}^{p}),
\end{equation}
where (recall that $K\leq1$)
$$
N_{1}=N_{1}(d,p,\theta,\tilde{\delta},\delta,\varepsilon,
\varepsilon_{1}),
\quad N_{2}=N_{2}(d,p,\theta,\tilde{\delta},\delta,\varepsilon,
\varepsilon_{1},\chi).
$$
Now we fix a $\chi=\chi(d,p,\theta,\tilde{\delta},\delta,
\varepsilon,\varepsilon_{1})
>0$ so that $N_{1}\chi\leq1/4$ and then find
a $\beta_{2}=\beta_{2}(d,p,\theta,\tilde{\delta},\delta,
\varepsilon,\varepsilon_{1})$
such that $N_{1}R\beta_{2}\leq1/4$. Then 
 estimate \eqref{6.2.7} will implies
\eqref{1.30.3} which along with Theorem \ref{theorem 1.30.1}
brings the proof of Theorem \ref{theorem 6.2.1}
to an end.

\end{document}